\documentclass{amsproc}%
\usepackage{graphicx}
\usepackage{amscd}
\usepackage{amsmath}
\usepackage{amsfonts}
\usepackage{amssymb}%
\theoremstyle{plain}

\numberwithin{equation}{section}

\begin{document}
\title[Centralizers in endomorphism rings]{Centralizers in endomorphism rings }
\author{Vesselin Drensky}
\address{Institute of Mathematics and Informatics, Bulgarian Academy of Sciences, 1113
Sofia, Bulgaria}
\email{drensky@math.bas.bg}
\author{Jen\H{o} Szigeti}
\address{Institute of Mathematics, University of Miskolc, Miskolc, Hungary 3515}
\email{jeno.szigeti@uni-miskolc.hu}
\author{Leon van Wyk}
\address{Department of Mathematical Sciences, Stellenbosch University\\
P/Bag X1, Matieland 7602, Stellenbosch, South Africa }
\email{lvw@sun.ac.za}
\thanks{The research of the first author was partially supported by Grant MI-1503/2005
of the Bulgarian National Science Fund. \noindent The second author was
supported by OTKA of Hungary No. K61007. The third author was supported by the
National Research Foundation of South Africa under Grant No. UID 61857. Any
opinion, findings and conclusions or recommendations expressed in this
material are those of the authors and therefore the National Research
Foundation does not accept any liability in regard thereto.}

\begin{abstract}
We prove that the centralizer Cen$(\varphi)\subseteq$ Hom$_{R}(M,M)$ of a
nilpotent endomorphism $\varphi$\ of a finitely generated semisimple left
$R$-module $_{R}M$ (over an arbitrary ring $R$) is the homomorphic image of
the opposite of a certain $Z(R)$-subalgebra of the full $m\times m$ matrix
algebra $M_{m}(R[z])$, where $m$ is the dimension (composition length) of
$\ker(\varphi)$. If $R$ is a local ring, then we provide an explicit
description of the above Cen$(\varphi)$. If in addition $Z(R)$ is a field and
$R/J(R)$ is finite dimensional over $Z(R)$, then we give a formula for the
$Z(R)$-dimension of Cen$(\varphi)$. If $R$ is a local ring, $\varphi$ is as
above and $\sigma\in$ Hom$_{R}(M,M)$ is arbitrary, then we give a complete
description of the containment Cen$(\varphi)\subseteq$ Cen$(\sigma)$ in terms
of an appropriate $R$-generating set of $_{R}M$. Using our results about
nilpotent endomorphisms, for an arbitrary (not necessarily nilpotent) linear
map $\varphi\in$ Hom$_{K}(V,V)$ of a finite dimensional vector space $V$\ over
a field $K$ we determine the PI-degree of Cen$(\varphi)$ and give other
information about the polynomial identities of Cen$(\varphi)$.

\end{abstract}
\maketitle

\noindent1. INTRODUCTION

\bigskip

Our work was motivated by one of the classical subjects of advanced linear
algebra. A detailed study of commuting matrices can be found in many of the
text books on linear algebra ([9, 10]). Commuting pairs and $k$-tuples of
$n\times n$\ matrices have been continuously in the focus of research (see,
for example, [5, 6]). If we replace $n\times n$ matrices by endomorphisms of
an $n$-generated module, we get a more general situation. The aim of the
present paper is to investigate the size and the PI properties of the
centralizer Cen$(\varphi)$\ of an element $\varphi$ in the endomorphism ring
Hom$_{R}(M,M)$.

In general, if $S$ is a ring (or algebra), then the centralizer%
\[
\text{Cen}(s)=\{u\in S\mid us=su\}
\]
of an element $s\in S$ is a subring (subalgebra) of $S$. Clearly, Cen$(s)$
satisfies all polynomial identities of $S$. However, Cen$(s)$ may also satisfy
some other polynomial identities. Thus the study of the PI properties of
Cen$(s)$, particularly in the case of classical rings (e.g. matrix rings),
deserves special attention.

Following the observations of Sections 2 and 3 about the nilpotent Jordan
normal base, Section 4 contains the following results about the centralizer
Cen$(\varphi)$ of a nilpotent endomorphism $\varphi$ of a finitely generated
semisimple left $R$-module $_{R}M$ (over an arbitrary ring $R$).

In Theorem 4.8 we prove that Cen$(\varphi)$ is the homomorphic image of the
opposite of a certain $Z(R)$-subalgebra of the full $m\times m$ matrix algebra
$M_{m}(R[z])$\ over the polynomial ring $R[z]$, where $m$ is the dimension
(composition length) of $\ker(\varphi)$. As a consequence, we obtain that
Cen$(\varphi)$\ satisfies all polynomial identities of $M_{m}^{\text{op}%
}(R[z])$; in particular, if $R$ is commutative, then the standard identity
$S_{2m}=0$ of degree $2m$ holds in Cen$(\varphi)$.

In Theorem 4.10 we obtain a complete characterization of Cen$(\varphi)$ if $R$
is a local ring.

In Theorem 4.11 we determine the $Z(R)$-dimension of Cen$(\varphi)$ in the
case when $R$ is local, $Z(R)$ is a field and $R/J(R)$ is finite dimensional
over $Z(R)$.

Section 5 contains further results about centralizers.

In Theorem 5.1 a complete description of Cen$(\varphi)$ is provided in terms
of an appropriate $R$-generating set of the finitely generated semisimple left
$R$-module $_{R}M$ (over an arbitrary ring $R$), with $\varphi
:M\longrightarrow M$ a so-called \textit{indecomposable nilpotent}
$R$-endomorphism. In particular, if $R$ is commutative, then we prove in
Corollary 5.2 that $\psi\in$ Cen $(\varphi)$ if and only if $\psi$ is a
polynomial expression of $\varphi$. A nilpotent linear map (or equivalently, a
nilpotent $n\times n$ matrix) is indecomposable if and only if its
characteristic and minimal polynomials coincide (see part (3) in Proposition
2.3). The following is a classical result about the centralizer (see [10]). If
$K$\ is a field and the characteristic and minimal polynomials of a (not
necessarily nilpotent) matrix $A\in M_{n}(K)$ coincide, then%
\[
\text{Cen}(A)=\{f(A)\mid f(z)\in K[z]\}.
\]
For an indecomposable nilpotent $\varphi$\ the mentioned corollary is a
generalization of the above result. We note that the above result on
Cen$(A)$\ is similar to Bergman's Theorem ([1]) about the centralizer of a
non-constant polynomial in the free associative algebra.

In Theorem 5.3 the containment relation Cen$(\varphi)\subseteq$ Cen$(\sigma)$
is considered. If $R$ is a local ring, $\varphi$ is nilpotent and $\sigma\in$
Hom$_{R}(M,M)$ is arbitrary, then we provide a complete description of the
situation Cen$(\varphi)\subseteq$ Cen$(\sigma)$ in terms of an appropriate
$R$-generating set of $_{R}M$. For two (not necessarily nilpotent) matrices
$A,B\in M_{n}(K)$, over an algebraically closed field $K$, the containment
Cen$(A)\subseteq$ Cen$(B)$ holds if and only if $B=f(A)$ for some $f(z)\in
K[z]$ (see [9] part VII, section 39). For a nilpotent $\varphi$\ our
description of Cen$(\varphi)\subseteq$ Cen$(\sigma)$ is a generalization of
the above result.

Section 6 is devoted to a study of the polynomial identities of the
centralizer. An arbitrary (not necessarily nilpotent) linear map of a finite
dimensional vector space, or equivalently, a matrix $A\in M_{n}(K)$ over a
field $K$ will be considered. First we show that it suffices to deal with
nilpotent matrices.

In Theorem 6.1 the radical and the semisimple component of the centralizer
Cen$(A)$ of a nilpotent matrix $A$\ is determined. The proof of Theorem 6.1 is
based on the explicit presentation of Cen$(A)$ in Theorem 4.10.

Corollary 6.2 is about the polynomial identities of the centralizer of a (not
necessarily nilpotent) matrix $A\in M_{n}(K)$. The Jordan normal form of $A$
over the algebraic closure of $K$ is considered. If $p$ is the maximum of the
number of elementary Jordan matrices of the same size and with the same
eigenvalue, and $T(S)\subseteq K\langle x_{1},x_{2},\ldots\rangle$ is the
T-ideal of the polynomial identities of the algebra $S$, then $T($%
Cen$(A))\supseteq T(M_{p}(K))^{q}$ for a suitable $q$, which also can be found
explicitly. Hence the PI-degree of Cen$(A)$ is equal to $p$.

Since all known results about matrix centralizers are closely connected with
the Jordan normal form, it is not surprising that our development depends on
the existence of the nilpotent Jordan normal base of a semisimple module with
respect to a given nilpotent endomorphism (guaranteed by one of the main
theorems of [11]).

\bigskip

\noindent2. THE NILPOTENT\ JORDAN NORMAL\ BASE

\bigskip

\noindent Throughout the paper a ring $R$\ means a (not necessarily
commutative) ring with identity, and $Z(R)$ and $J(R)$ denote the centre and
the Jacobson radical of $R$, respectively. Let $\varphi:M\longrightarrow M$ be
an $R$-endomorphism of the (unitary) left $R$-module $_{R}M$. A subset%
\[
\{x_{\gamma,i}\mid\gamma\in\Gamma,1\leq i\leq k_{\gamma}\}\subseteq M
\]
is called a \textit{nilpotent Jordan normal base} of $_{R}M$\ with respect to
$\varphi$ if each $R$-submodule $Rx_{\gamma,i}\leq M$ is simple,%
\[
\underset{\gamma\in\Gamma,1\leq i\leq k_{\gamma}}{\oplus}Rx_{\gamma,i}=M
\]
is a direct sum, $\varphi(x_{\gamma,i})=x_{\gamma,i+1}$ $,$ $\varphi
(x_{\gamma,k_{\gamma}})=x_{\gamma,k_{\gamma}+1}=0$ for all\textit{ }$\gamma
\in\Gamma$, $1\leq i\leq k_{\gamma}$, and the set $\{k_{\gamma}\mid\gamma
\in\Gamma\}$ of integers is bounded. For $i\geq k_{\gamma}+1$ we assume that
$x_{\gamma,i}=0$ holds in $M$. Now $\Gamma$\ is called the set of (Jordan-)
blocks and the size of the block $\gamma\in\Gamma$ is the integer $k_{\gamma
}\geq1$. Obviously, the existence of a nilpotent Jordan normal base implies
that $_{R}M$ is semisimple and $\varphi$ is nilpotent with $\varphi^{n}%
=0\neq\varphi^{n-1}$, where%
\[
n=\max\{k_{\gamma}\mid\gamma\in\Gamma\}
\]
is the index of nilpotency. If $_{R}M$ is finitely generated, then%
\[
\underset{\gamma\in\Gamma}{\sum}k_{\gamma}=\dim_{R}(M)
\]
is the dimension of $_{R}M$\ (equivalently: the composition length of $_{R}M$
or the height of the submodule lattice of $_{R}M$). Clearly,%
\[
\varphi(M)=\varphi\left(  \underset{\gamma\in\Gamma,1\leq i\leq k_{\gamma}%
}{\sum}Rx_{\gamma,i}\right)  =\underset{\gamma\in\Gamma,1\leq i\leq k_{\gamma
}}{\sum}R\varphi(x_{\gamma,i})=\underset{\gamma\in\Gamma,1\leq i\leq
k_{\gamma}-1}{\sum}Rx_{\gamma,i+1}%
\]
implies that%
\[
\text{im}(\varphi)=\underset{\gamma\in\Gamma,1\leq i\leq k_{\gamma}-1}{\oplus
}Rx_{\gamma,i+1}=\underset{\gamma\in\Gamma^{\prime},2\leq i^{\prime}\leq
k_{\gamma}}{\oplus}Rx_{\gamma,i^{\prime}},
\]
where%
\[
\Gamma^{\prime}=\{\gamma\in\Gamma\mid k_{\gamma}\geq2\}\text{ and }%
\Gamma\setminus\Gamma^{\prime}=\{\gamma\in\Gamma\mid k_{\gamma}=1\}.
\]
Any element $u\in M$ can be written as%
\[
u=\underset{\gamma\in\Gamma,1\leq i\leq k_{\gamma}}{\sum}a_{\gamma,i}%
x_{\gamma,i},
\]
where $\{(\gamma,i)\mid\gamma\in\Gamma,1\leq i\leq k_{\gamma},$ and
$a_{\gamma,i}\neq0\}$ is finite and all summands $a_{\gamma,i}x_{\gamma,i}$
are uniquely determined by $u$. Since%
\[
\varphi(u)=\underset{\gamma\in\Gamma,1\leq i\leq k_{\gamma}}{\sum}a_{\gamma
,i}\varphi(x_{\gamma,i})=\underset{\gamma\in\Gamma,1\leq i\leq k_{\gamma}%
}{\sum}a_{\gamma,i}x_{\gamma,i+1}=0
\]
is equivalent to the condition that $a_{\gamma,i}x_{\gamma,i+1}=0$ for all
$\gamma\in\Gamma,1\leq i\leq k_{\gamma}-1$, we obtain that%
\[
\varphi(u)=0\Longleftrightarrow u=\underset{\gamma\in\Gamma}{\sum}%
a_{\gamma,k_{\gamma}}x_{\gamma,k_{\gamma}}.
\]
Indeed, $a_{\gamma,i}x_{\gamma,i}\neq0$ ($1\leq i\leq k_{\gamma}-1$) would
imply that $ba_{\gamma,i}x_{\gamma,i}=x_{\gamma,i}$ for some $b\in R$ (note
that $Rx_{\gamma,i}\leq M$ is simple), whence%
\[
x_{\gamma,i+1}=\varphi(x_{\gamma,i})=ba_{\gamma,i}\varphi(x_{\gamma
,i})=ba_{\gamma,i}x_{\gamma,i+1}=0
\]
can be derived, a contradiction. It follows that%
\[
\ker(\varphi)=\underset{\gamma\in\Gamma}{\oplus}Rx_{\gamma,k_{\gamma}}%
\]
and $\dim_{R}(\ker(\varphi))=\left\vert \Gamma\right\vert $ in case of a
finite $\Gamma$. The following is one of the main results in [11].

\bigskip

\noindent\textbf{2.1.Theorem.}\textit{ Let }$\varphi:M\longrightarrow
M$\textit{ be an }$R$\textit{-endomorphism of the left }$R$\textit{-module
}$_{R}M$\textit{. Then the following are equivalent.}

\begin{enumerate}
\item $_{R}M$\textit{ is a semisimple left }$R$\textit{-module and }$\varphi
$\textit{\ is nilpotent.}

\item \textit{There exists a nilpotent Jordan normal base of }$_{R}%
M$\textit{\ with respect to }$\varphi$\textit{.}
\end{enumerate}

\bigskip

\noindent\textbf{2.2.Proposition.}\textit{ Let }$\varphi:M\longrightarrow
M$\textit{ be a nilpotent }$R$\textit{-endomorphism of the finitely generated
semisimple left }$R$\textit{-module }$_{R}M$\textit{. If }$\{x_{\gamma,i}%
\mid\gamma\in\Gamma,1\leq i\leq k_{\gamma}\}$\textit{ and }$\{y_{\delta,j}%
\mid\delta\in\Delta,1\leq j\leq l_{\delta}\}$\textit{\ are nilpotent Jordan
normal bases of }$_{R}M$\textit{\ with respect to }$\varphi$\textit{, then
there exists a bijection }$\pi:\Gamma\longrightarrow\Delta$\textit{ such that
}$k_{\gamma}=l_{\pi(\gamma)}$\textit{ for all }$\gamma\in\Gamma$\textit{. Thus
the sizes of the blocks of a nilpotent Jordan normal base are unique up to a
permutation of the blocks.}

\bigskip

\noindent\textbf{Proof.} We apply induction on the index of the nilpotency of
$\varphi$. If $\varphi=0$, then we have $k_{\gamma}=l_{\delta}=1$ for all
$\gamma\in\Gamma$, $\delta\in\Delta$, and%
\[
\underset{\gamma\in\Gamma}{\oplus}Rx_{\gamma,1}=\underset{\delta\in\Delta
}{\oplus}Ry_{\delta,1}=M
\]
implies the existence of a bijection $\pi:\Gamma\longrightarrow\Delta$
(Krull-Schmidt, Kurosh-Ore). Assume that our statement holds for any
$R$-endomorphism $\phi:N\longrightarrow N$ with $_{R}N$ being a finitely
generated semisimple left $R$-module and $\phi^{n-2}\neq0=\phi^{n-1}$.
Consider the situation described in the proposition with $\varphi^{n-1}%
\neq0=\varphi^{n}$, then%
\[
\text{im}(\varphi)=\underset{\gamma\in\Gamma^{\prime},2\leq i^{\prime}\leq
k_{\gamma}}{\oplus}Rx_{\gamma,i^{\prime}}%
\]
ensures that%
\[
\{x_{\gamma,i^{\prime}}\mid\gamma\in\Gamma^{\prime},2\leq i^{\prime}\leq
k_{\gamma}\}
\]
is a nilpotent Jordan normal base of the left $R$-submodule im$(\varphi)\leq
M$ of $_{R}M$\ with respect to the restricted $R$-endomorphism $\varphi:$
im$(\varphi)\longrightarrow$ im$(\varphi)$. The same holds for%
\[
\{y_{\delta,j^{\prime}}\mid\delta\in\Delta^{\prime},2\leq j^{\prime}\leq
l_{\delta}\}.
\]
Since we have $\phi^{n-2}\neq0=\phi^{n-1}$ for $\phi=\varphi\upharpoonright$
im$(\varphi)$, our assumption ensures the existence of a bijection $\pi
:\Gamma^{\prime}\longrightarrow\Delta^{\prime}$ such that $k_{\gamma}%
-1=l_{\pi(\gamma)}-1$ for all $\gamma\in\Gamma^{\prime}$. In view of%
\[
\ker(\varphi)=\underset{\gamma\in\Gamma}{\oplus}Rx_{\gamma,k_{\gamma}%
}=\underset{\delta\in\Delta}{\oplus}Ry_{\delta,l_{\delta}}%
\]
we obtain that $\left\vert \Gamma\right\vert =\left\vert \Delta\right\vert $
(Krull-Schmidt, Kurosh-Ore), whence $\left\vert \Gamma\setminus\Gamma^{\prime
}\right\vert =\left\vert \Delta\setminus\Delta^{\prime}\right\vert $ follows.
Thus we have a bijection $\pi^{\ast}:\Gamma\setminus\Gamma^{\prime
}\longrightarrow\Delta\setminus\Delta^{\prime}$ and the natural map%
\[
\pi\sqcup\pi^{\ast}:\Gamma^{\prime}\cup(\Gamma\setminus\Gamma^{\prime
})\longrightarrow\Delta^{\prime}\cup(\Delta\setminus\Delta^{\prime})
\]
is a bijection with the desired property. $\square$

\bigskip

\noindent We call a nilpotent element $s\in S$ of the ring $S$%
\ \textit{decomposable} if $es=se$ holds for some idempotent element $e\in S$
($e^{2}=e$) with $0\neq e\neq1$. A nilpotent element which is not decomposable
is called \textit{indecomposable}.

\bigskip

\noindent\textbf{2.3.Proposition.}\textit{ Let }$\varphi:M\longrightarrow
M$\textit{ be a nonzero nilpotent }$R$\textit{-endomorphism of the semisimple
left }$R$\textit{-module }$_{R}M$\textit{. Then the following are equivalent.}

\begin{enumerate}
\item \textit{There is a nilpotent Jordan normal base }$\{x_{i}\mid1\leq i\leq
n\}$\textit{\ of }$_{R}M$\textit{\ with respect to }$\varphi$\textit{
consisting of one block (thus }$\left\vert \Gamma\right\vert =1$\textit{\ for
any nilpotent Jordan normal base }$\{x_{\gamma,i}\mid\gamma\in\Gamma,1\leq
i\leq k_{\gamma}\}$\textit{\ of }$_{R}M$\textit{\ with respect to }$\varphi
$\textit{).}

\item $\varphi$\textit{ is an indecomposable nilpotent element of the ring
}Hom$_{R}(M,M)$\textit{.}

\item $_{R}M$\textit{ is finitely generated and }$\varphi^{d-1}\neq0$\textit{,
where }$d=\dim_{R}(M)$\textit{ is the dimension of }$_{R}M$.
\end{enumerate}

\bigskip

\noindent\textbf{Proof.}

\noindent(1)$\Longrightarrow$(3): Clearly,%
\[
\underset{1\leq i\leq n}{\oplus}Rx_{i}=M
\]
implies that we have $d=n$\ for the dimension of $_{R}M$, whence%
\[
\varphi^{d-1}(x_{1})=\varphi^{n-1}(x_{1})=x_{n}\neq0
\]
follows.

\noindent(3)$\Longrightarrow$(1): Let $\{x_{\gamma,i}\mid\gamma\in\Gamma,1\leq
i\leq k_{\gamma}\}$ be a nilpotent Jordan normal base of $_{R}M$\ with respect
to $\varphi$. Suppose that $\left\vert \Gamma\right\vert \geq2$, then%
\[
n=\max\{k_{\gamma}\mid\gamma\in\Gamma\}\leq d-1,
\]
where $d=\underset{\gamma\in\Gamma}{\sum}k_{\gamma}=\dim_{R}(M)$. Thus
$\varphi^{d-1}=\varphi^{(d-1)-n}\circ\varphi^{n}=0$, a contradiction.

\noindent(1)$\Longrightarrow$(2): Suppose that $\varepsilon\circ
\varphi=\varphi\circ\varepsilon$ holds for some idempotent endomorphism
$\varepsilon\in$ Hom$_{R}(M,M)$ with $0\neq\varepsilon\neq1$. Then%
\[
\text{im}(\varepsilon)\oplus\text{im}(1-\varepsilon)=M
\]
for the non-zero (semisimple) $R$-submodules im$(\varepsilon)$ and
im$(1-\varepsilon)$ of $_{R}M$. Now $\varepsilon\circ\varphi=\varphi
\circ\varepsilon$ ensures that $\varphi:$ im$(\varepsilon)\longrightarrow$
im$(\varepsilon)$ and $\varphi:$ im$(1-\varepsilon)\longrightarrow$
im$(1-\varepsilon)$. Since these restricted $R$-endomorphisms are nilpotent,
we have a nilpotent Jordan normal base of im$(\varepsilon)$ with respect to
$\varphi\upharpoonright$ im$(\varepsilon)$ and a nilpotent Jordan normal base
of im$(1-\varepsilon)$ with respect to $\varphi\upharpoonright$
im$(1-\varepsilon)$. The union of these two bases gives a nilpotent Jordan
normal base of $M$ with respect to $\varphi$ consisting of more than one
block, a contradiction (the direct sum property of the new base is a
consequence of the modularity of the submodule lattice of $_{R}M$).

\noindent(2)$\Longrightarrow$(1): Suppose that $\{x_{\gamma,i}\mid\gamma
\in\Gamma,1\leq i\leq k_{\gamma}\}$ is a nilpotent Jordan normal base of
$_{R}M$\ with respect to $\varphi$ with $\left\vert \Gamma\right\vert \geq2$
and fix an element $\delta\in\Gamma$. Consider the non-zero $\varphi
$-invariant $R$-submodules
\[
N_{\delta}^{\prime}=\underset{1\leq i\leq k_{\delta}}{\oplus}Rx_{\delta
,i}\text{ and }N_{\delta}^{\prime\prime}=\underset{\gamma\in\Gamma
\setminus\{\delta\},1\leq i\leq k_{\gamma}}{\oplus}Rx_{\gamma,i}\text{ },
\]
then $M=N_{\delta}^{\prime}\oplus N_{\delta}^{\prime\prime}$ and define
$\varepsilon_{\delta}:M\longrightarrow M$ as the natural projection of $M$
onto $N_{\delta}^{\prime}$. Then $\varepsilon_{\delta}(u)=u^{\prime}$, where
$u=u^{\prime}+u^{\prime\prime}$ is the unique sum presentation of $u\in M$
with $u^{\prime}\in N_{\delta}^{\prime}$ and $u^{\prime\prime}\in N_{\delta
}^{\prime\prime}$. It is straightforward to see that $\varepsilon_{\delta
}\circ\varepsilon_{\delta}=\varepsilon_{\delta}$, $0\neq\varepsilon_{\delta
}\neq1$\ and $\varepsilon_{\delta}\circ\varphi=\varphi\circ\varepsilon
_{\delta}$ hold. $\square$

\bigskip

\noindent3. THE\ MODULE\ STRUCTURE\ INDUCED\ BY\ AN ENDOMORPHISM

\bigskip

\noindent Let $R[z]$ denote the ring of polynomials of the commuting
indeterminate $z$ with coefficients in $R$. The ideal $(z^{k})=R[z]z^{k}%
=z^{k}R[z]\vartriangleleft R[z]$ generated by $z^{k}$ will be considered in
the sequel. If $\varphi:M\longrightarrow M$ is an arbitrary $R$-endomorphism
of the left $R$-module $_{R}M$, then for $u\in M$ and%
\[
f(z)=a_{1}+a_{2}z+\cdots+a_{n+1}z^{n}\in R[z]
\]
(unusual use of indices!) the left multiplication%
\[
f(z)\ast u=a_{1}u+a_{2}\varphi(u)+\cdots+a_{n+1}\varphi^{n}(u)
\]
defines a natural left $R[z]$-module structure on $M$. This left action of
$R[z]$ on $M$ extends the left action of $R$ on $_{R}M$. Note that%
\[
z^{n}\ast u=\varphi^{n}(u)\text{ and }\varphi(f(z)\ast u)=(zf(z))\ast u.
\]
For any $R$-endomorphism $\psi\in$ Hom$_{R}(M,M)$ with $\psi\circ
\varphi=\varphi\circ\psi$ we have%
\[
\psi(f(z)\ast u)=f(z)\ast\psi(u)
\]
and hence $\psi:M\longrightarrow M$ is an $R[z]$-endomorphism of the left
$R[z]$-module $_{R[z]}M$. On the other hand, if $\psi:M\longrightarrow M$ is
an $R[z]$-endomorphism of $_{R[z]}M$, then%
\[
\psi(\varphi(u))=\psi(z\ast u)=z\ast\psi(u)=\varphi(\psi(u))
\]
implies that $\psi\circ\varphi=\varphi\circ\psi$. The centralizer
\[
\text{Cen}(\varphi)=\{\psi\mid\psi\in\text{Hom}_{R}(M,M)\text{ and }\psi
\circ\varphi=\varphi\circ\psi\}
\]
of $\varphi$ is a $Z(R)$-subalgebra of Hom$_{R}(M,M)$ and the argument above
gives that%
\[
\text{Cen}(\varphi)=\text{Hom}_{R[z]}(M,M).
\]

\noindent For a set $\Gamma\neq\varnothing$, the $\Gamma$-copower
$\underset{\gamma\in\Gamma}{%
{\displaystyle\coprod}
}R[z]$\ of the ring $R[z]$ is an ideal of the $\Gamma$-direct power ring
$\underset{\gamma\in\Gamma}{%
{\displaystyle\prod}
}R[z]$ consisting of all elements $\mathbf{f}=(f_{\gamma}(z))_{\gamma\in
\Gamma}$ with a finite set $\{\gamma\in\Gamma\mid f_{\gamma}(z)\neq0\}$\ of
non-zero coordinates. The power (copower) has a natural $(R[z],R[z])$-bimodule
structure. If $\Gamma$ is finite, then%
\[
(R[z])^{\Gamma}=\underset{\gamma\in\Gamma}{%
{\displaystyle\coprod}
}R[z]=\underset{\gamma\in\Gamma}{%
{\displaystyle\prod}
}R[z].
\]
\noindent If $\{x_{\gamma,i}\mid\gamma\in\Gamma,1\leq i\leq k_{\gamma}\}$ is a
nilpotent Jordan normal base of $_{R}M$\ with respect to a nilpotent
endomorphism $\varphi$, then for an element $\mathbf{f}=(f_{\gamma
}(z))_{\gamma\in\Gamma}$ with%
\[
f_{\gamma}(z)=a_{\gamma,1}+a_{\gamma,2}z+\cdots+a_{\gamma,n_{\gamma}%
+1}z^{n_{\gamma}}%
\]
the formula%
\[
\Phi(\mathbf{f})=\underset{\gamma\in\Gamma,1\leq i\leq k_{\gamma}}{%
{\displaystyle\sum}
}a_{\gamma,i}x_{\gamma,i}=\underset{\gamma\in\Gamma}{%
{\displaystyle\sum}
}f_{\gamma}(z)\ast x_{\gamma,1}%
\]
defines a function%
\[
\Phi:\underset{\gamma\in\Gamma}{%
{\displaystyle\coprod}
}R[z]\longrightarrow M.
\]

\bigskip

\noindent\textbf{3.1.Proposition.}\textit{ For a nilpotent endomorphism
}$\varphi\in$ Hom$_{R}(M,M)$\textit{\ of the semisimple left }$R$%
\textit{-module }$_{R}M$\textit{, the function }$\Phi$\textit{\ is a
surjective left }$R[z]$\textit{-homomorphism. We have }$\varphi(\Phi
(\mathbf{f}))=\Phi(z\mathbf{f})$\textit{ for all }$\mathbf{f}\in
\underset{\gamma\in\Gamma}{%
{\displaystyle\coprod}
}R[z]$\textit{ and the kernel}%
\[
\underset{\gamma\in\Gamma}{%
{\displaystyle\coprod}
}J(R)[z]+(z^{k_{\gamma}})\subseteq\ker(\Phi)\vartriangleleft_{l}%
\underset{\gamma\in\Gamma}{%
{\displaystyle\prod}
}R[z]
\]
\textit{is a left ideal of the power (and hence of the copower) ring. If }%
$R$\textit{ is a local ring (}$R/J(R)$\textit{ is a division ring), then}%
\[
\underset{\gamma\in\Gamma}{%
{\displaystyle\coprod}
}J(R)[z]+(z^{k_{\gamma}})=\ker(\Phi).
\]

\bigskip

\noindent\textbf{Proof.} Clearly,%
\[
\underset{\gamma\in\Gamma,1\leq i\leq k_{\gamma}}{\sum}Rx_{\gamma,i}=M
\]
implies that $\Phi$ is surjective. The second part of the defining formula
gives that $\Phi$ is a left $R[z]$-homomorphism:%
\[
\Phi(g(z)\mathbf{f})=\underset{\gamma\in\Gamma}{%
{\displaystyle\sum}
}(g(z)f_{\gamma}(z))\ast x_{\gamma,1}=\underset{\gamma\in\Gamma}{%
{\displaystyle\sum}
}g(z)\ast(f_{\gamma}(z)\ast x_{\gamma,1})=g(z)\ast\Phi(\mathbf{f}),
\]
where $g(z)\in R[z]$. We also have%
\[
\varphi(\Phi(\mathbf{f}))=\underset{\gamma\in\Gamma}{%
{\displaystyle\sum}
}\varphi(f_{\gamma}(z)\ast x_{\gamma,1})=\underset{\gamma\in\Gamma}{%
{\displaystyle\sum}
}(zf_{\gamma}(z))\ast x_{\gamma,1}=\Phi(z\mathbf{f}).
\]
If $\mathbf{f}\in\underset{\gamma\in\Gamma}{%
{\displaystyle\coprod}
}J(R)[z]+(z^{k_{\gamma}})$, then%
\[
f_{\gamma}(z)=(a_{\gamma,1}+a_{\gamma,2}z+\cdots+a_{\gamma,k_{\gamma}%
}z^{k_{\gamma}-1})+(a_{\gamma,k_{\gamma}+1}z^{k_{\gamma}}+\cdots
+a_{\gamma,n_{\gamma}+1}z^{n_{\gamma}})
\]
with $a_{\gamma,i}\in J(R)$, $1\leq i\leq k_{\gamma}$. Since $Rx_{\gamma,i}$
is simple, we have $J(R)x_{\gamma,i}=\{0\}$. Thus $\varphi^{k_{\gamma}%
}(x_{\gamma,1})=0$ implies that $f_{\gamma}(z)\ast x_{\gamma,1}=0$, whence
$\Phi(\mathbf{f})=0$ follows. Take an element $\mathbf{g}=(g_{\gamma
}(z))_{\gamma\in\Gamma}$ of the direct power and suppose that $\Phi
(\mathbf{f})=0$ in $M$. Then%
\[
\underset{\gamma\in\Gamma,1\leq i\leq k_{\gamma}}{\oplus}Rx_{\gamma,i}=M
\]
implies that $a_{\gamma,i}x_{\gamma,i}=0$ for all $\gamma\in\Gamma$ and $1\leq
i\leq k_{\gamma}$. Thus $f_{\gamma}(z)\ast x_{\gamma,1}=0$ for all $\gamma
\in\Gamma$. It follows that%
\[
\Phi(\mathbf{gf})=\underset{\gamma\in\Gamma}{%
{\displaystyle\sum}
}(g_{\gamma}(z)f_{\gamma}(z))\ast x_{\gamma,1}=\underset{\gamma\in\Gamma}{%
{\displaystyle\sum}
}g_{\gamma}(z)\ast(f_{\gamma}(z)\ast x_{\gamma,1})=0,
\]
whence $\mathbf{gf}\in\ker(\Phi)$ can be deduced.

\noindent If $R$ is a local ring and $a_{\gamma,i}x_{\gamma,i}=0$ for some
$1\leq i\leq k_{\gamma}$, then $a_{\gamma,i}\in J(R)$. Thus $\Phi
(\mathbf{f})=0$ implies that%
\[
f_{\gamma}(z)=(a_{\gamma,1}+a_{\gamma,2}z+\cdots+a_{\gamma,k_{\gamma}%
}z^{k_{\gamma}-1})+(a_{\gamma,k_{\gamma}+1}z^{k_{\gamma}}+\cdots
+a_{\gamma,n_{\gamma}+1}z^{n_{\gamma}})\in J(R)[z]+(z^{k_{\gamma}}).
\]
It follows that%
\[
\mathbf{f}\in\underset{\gamma\in\Gamma}{%
{\displaystyle\coprod}
}J(R)[z]+(z^{k_{\gamma}}).\ \square
\]

\bigskip

\noindent4. THE\ CENTRALIZER\ OF\ A\ NILPOTENT ENDOMORPHISM

\bigskip

\noindent Let $X=\{x_{\gamma,i}\mid\gamma\in\Gamma,1\leq i\leq k_{\gamma}\}$
be a nilpotent Jordan normal base of $_{R}M$\ with respect to the nilpotent
endomorphism $\varphi\in$ Hom$_{R}(M,M)$. We keep the notations of the
previous section and in the rest of the paper we assume that $_{R}M$ is
finitely generated, i.e. that $\Gamma$\ is finite.

\noindent A linear order on $\Gamma$, say $\Gamma=\{1,2,\ldots,m\}$, allows us
to view an element $\mathbf{f}=(f_{\gamma}(z))_{\gamma\in\Gamma}$ of
$(R[z])^{\Gamma}$ as a $1\times m$ matrix (a row vector) over $R[z]$. For an
$m\times m$ matrix $\mathbf{P}=[p_{\delta,\gamma}(z)]$ in $M_{m}(R[z])$ the
matrix product%
\[
\mathbf{fP=}\underset{\delta\in\Gamma}{%
{\displaystyle\sum}
}f_{\delta}(z)\mathbf{p}_{\delta}%
\]
of $\mathbf{f}$ and $\mathbf{P}$ is a $1\times m$ matrix (row vector) in
$(R[z])^{\Gamma}$, where $\mathbf{p}_{\delta}=(p_{\delta,\gamma}%
(z))_{\gamma\in\Gamma}$ is the $\delta$-th row vector of $\mathbf{P}$ and%
\[
(\mathbf{fP})_{\gamma}=\underset{\delta\in\Gamma}{%
{\displaystyle\sum}
}f_{\delta}(z)p_{\delta,\gamma}(z).
\]
We define the subsets%
\[
\mathcal{I}(X)=\{\mathbf{P}\in M_{m}(R[z])\mid\mathbf{P}=[p_{\delta,\gamma
}(z)]\text{ and }p_{\delta,\gamma}(z)\in J(R)[z]+(z^{k_{\gamma}})\text{ for
all }\delta,\gamma\in\Gamma\},
\]%
\[
\mathcal{N}(X)=\{\mathbf{P}\in M_{m}(R[z])\mid\mathbf{P}=[p_{\delta,\gamma
}(z)]\text{ and }z^{k_{\delta}}p_{\delta,\gamma}(z)\in J(R)[z]+(z^{k_{\gamma}%
})\text{ for all }\delta,\gamma\in\Gamma\}
\]
and%
\[
\mathcal{M}(X)=\{\mathbf{P}\in M_{m}(R[z])\mid\mathbf{fP}\in\ker(\Phi)\text{
for all }\mathbf{f}\in\ker(\Phi)\}
\]
of $M_{m}(R[z])$. Note that $\mathcal{I}(X)$ and $\mathcal{N}(X)$ are
$(R[z],R[z])$-sub-bimodules of $M_{m}(R[z])$ in a natural way.

\noindent For $\delta,\gamma\in\Gamma$ let $k_{\delta,\gamma}=k_{\gamma
}-k_{\delta}$ when $1\leq k_{\delta}<k_{\gamma}\leq n$ and $k_{\delta,\gamma
}=0$ otherwise. It can be verified, that the condition $z^{k_{\delta}%
}p_{\delta,\gamma}(z)\in J(R)[z]+(z^{k_{\gamma}})$ is equivalent to
$p_{\delta,\gamma}(z)\in J(R)[z]+(z^{k_{\delta,\gamma}})$. Note that $z^{0}=1$
and $(z^{0})=R[z]$.

\bigskip

\noindent\textbf{4.1.Remark.}\textit{ If }$\mathbf{E}_{\delta,\gamma}$\textit{
denotes the }$m\times m$\textit{\ standard matrix unit over }$R[z]$%
\textit{\ with }$1$\textit{ in the }$(\delta,\gamma)$\textit{ entry and zeros
in the other entries, then }$\mathbf{E}_{\delta,\gamma}\in\mathcal{N}%
(X)$\textit{ for all }$\delta,\gamma\in\Gamma$\textit{ with }$k_{\delta}\geq
k_{\gamma}$\textit{.}

\bigskip

\noindent\textbf{4.2.Lemma.}\textit{ }$\mathcal{I}(X)\vartriangleleft_{l}%
M_{m}(R[z])$\textit{ is a left ideal, }$\mathcal{N}(X)\subseteq M_{m}%
(R[z])$\textit{ is a subring, }$\mathcal{I}(X)\vartriangleleft\mathcal{N}%
(X)$\textit{ is an ideal and }$\mathcal{M}(X)$\textit{\ is a }$Z(R)$%
\textit{-subalgebra of }$M_{m}(R[z])$\textit{. If }$R$\textit{ is a local
ring, then }$\mathcal{N}(X)=\mathcal{M}(X)$\textit{.}

\bigskip

\noindent\textbf{Proof.} For the elements $\mathbf{P}=[p_{\delta,\gamma}(z)]$
and $\mathbf{Q}=[q_{\delta,\gamma}(z)]$ of $M_{m}(R[z])$ take $\mathbf{PQ}%
=[r_{\delta,\gamma}(z)]$, where%
\[
r_{\delta,\gamma}(z)=\underset{\tau\in\Gamma}{\sum}p_{\delta,\tau}%
(z)q_{\tau,\gamma}(z).
\]
If $\mathbf{Q}\in\mathcal{I}(X)$, then $q_{\tau,\gamma}(z)\in
J(R)[z]+(z^{k_{\gamma}})$ and $J(R)[z]+(z^{k_{\gamma}})\vartriangleleft R[z]$
imply that $r_{\delta,\gamma}(z)\in J(R)[z]+(z^{k_{\gamma}})$. Thus
$\mathbf{PQ}\in\mathcal{I}(X)$ and $\mathcal{I}(X)$ is a left ideal of
$M_{m}(R[z])$.

\noindent If $\mathbf{P,Q}\in\mathcal{N}(X)$, then $z^{k_{\delta}}%
p_{\delta,\tau}(z)\in J(R)[z]+(z^{k_{\tau}})$ and $z^{k_{\tau}}q_{\tau,\gamma
}(z)\in J(R)[z]+(z^{k_{\gamma}})$ for all $\delta,\tau,\gamma\in\Gamma$. It
follows that%
\[
z^{k_{\delta}}p_{\delta,\tau}(z)=g_{\delta,\tau}(z)+z^{k_{\tau}}h_{\delta
,\tau}(z)
\]
with $g_{\delta,\tau}(z)\in J(R)[z]$ and $h_{\delta,\tau}(z)\in R[z]$. Thus%
\[
z^{k_{\delta}}r_{\delta,\gamma}(z)=\underset{\tau\in\Gamma}{\sum}z^{k_{\delta
}}p_{\delta,\tau}(z)q_{\tau,\gamma}(z)=\underset{\tau\in\Gamma}{\sum
}(g_{\delta,\tau}(z)+z^{k_{\tau}}h_{\delta,\tau}(z))q_{\tau,\gamma}(z)=
\]%
\[
\underset{\tau\in\Gamma}{\sum}(g_{\delta,\tau}(z)q_{\tau,\gamma}%
(z)+h_{\delta,\tau}(z)z^{k_{\tau}}q_{\tau,\gamma}(z))
\]
and $z^{k_{\tau}}q_{\tau,\gamma}(z)\in J(R)[z]+(z^{k_{\gamma}})$ ensure that
$z^{k_{\delta}}r_{\delta,\gamma}(z)$ is in $J(R)[z]+(z^{k_{\gamma}})$.
Consequently, we obtain that $\mathbf{PQ}\in\mathcal{N}(X)$. Hence
$\mathcal{N}(X)$ is a subring of $M_{m}(R[z])$.

\noindent If $\mathbf{P}\in\mathcal{I}(X)$ and $\mathbf{Q}\in\mathcal{N}(X)$,
then $p_{\delta,\tau}(z)\in J(R)[z]+(z^{k_{\tau}})$ and $z^{k_{\tau}}%
q_{\tau,\gamma}(z)\in J(R)[z]+(z^{k_{\gamma}})$ for all $\delta,\tau,\gamma
\in\Gamma$. It follows that
\[
p_{\delta,\tau}(z)=u_{\delta,\tau}(z)+z^{k_{\tau}}v_{\delta,\tau}(z)
\]
with $u_{\delta,\tau}(z)\in J(R)[z]$ and $v_{\delta,\tau}(z)\in R[z]$. Thus%
\[
r_{\delta,\gamma}(z)=\underset{\tau\in\Gamma}{\sum}p_{\delta,\tau}%
(z)q_{\tau,\gamma}(z)=\underset{\tau\in\Gamma}{\sum}(u_{\delta,\tau
}(z)+z^{k_{\tau}}v_{\delta,\tau}(z))q_{\tau,\gamma}(z)=
\]%
\[
\underset{\tau\in\Gamma}{\sum}(u_{\delta,\tau}(z)q_{\tau,\gamma}%
(z)+v_{\delta,\tau}(z)z^{k_{\tau}}q_{\tau,\gamma}(z))
\]
and $z^{k_{\tau}}q_{\tau,\gamma}(z)\in J(R)[z]+(z^{k_{\gamma}})$ ensure that
$r_{\delta,\gamma}(z)$ is in $J(R)[z]+(z^{k_{\gamma}})$. Consequently, we
obtain that $\mathbf{PQ}\in\mathcal{I}(X)$. Hence $\mathcal{I}(X)$ is an ideal
of $\mathcal{N}(X)$.

\noindent If $\mathbf{P},\mathbf{Q}\in\mathcal{M}(X)$, $\mathbf{f}\in\ker
(\Phi)$ and $c\in Z(R)$, then%
\[
\Phi(\mathbf{f}(c\mathbf{P}))=\Phi(c(\mathbf{fP}))=c\Phi(\mathbf{fP})=0,
\]%
\[
\Phi(\mathbf{f}(\mathbf{P}+\mathbf{Q}))=\Phi(\mathbf{fP}+\mathbf{fQ}%
)=\Phi(\mathbf{fP})+\Phi(\mathbf{fQ})=0
\]
and $\mathbf{fP}\in\ker(\Phi)$ implies that%
\[
\Phi(\mathbf{f}(\mathbf{PQ}))=\Phi((\mathbf{fP})\mathbf{Q})=0,
\]
whence $\mathbf{f}(\mathbf{PQ})\in\ker(\Phi)$ follows. Thus $c\mathbf{P}%
,\mathbf{P}+\mathbf{Q},\mathbf{PQ}\in\mathcal{M}(X)$, proving that
$\mathcal{M}(X)$\ is a $Z(R)$-subalgebra of $M_{m}(R[z])$.

\noindent If $R$\ is a local ring, then Proposition 3.1 gives that%
\[
\ker(\Phi)=\underset{\gamma\in\Gamma}{%
{\displaystyle\coprod}
}J(R)[z]+(z^{k_{\gamma}}).
\]
Now $\mathbf{e}_{\delta}\in\ker(\Phi)$, where $\mathbf{e}_{\delta}$ denotes
the vector with $z^{k_{\delta}}$ in its $\delta$-coordinate and zeros in all
other places.

\noindent If $\mathbf{P}\in\mathcal{M}(X)$, then $\mathbf{e}_{\delta
}\mathbf{P}\in\ker(\Phi)$ implies that $z^{k_{\delta}}p_{\delta,\gamma}(z)\in
J(R)[z]+(z^{k_{\gamma}})$ for all $\delta,\gamma\in\Gamma$, whence
$\mathbf{P}\in\mathcal{N}(X)$ follows.

\noindent If $\mathbf{P}\in\mathcal{N}(X)$ and $\mathbf{f}\in\ker(\Phi)$, then
$z^{k_{\delta}}p_{\delta,\gamma}(z)\in J(R)[z]+(z^{k_{\gamma}})$ for all
$\delta,\gamma\in\Gamma$ and $f_{\gamma}(z)=g_{\gamma}(z)+z^{k_{\gamma}%
}h_{\gamma}(z)$ with $g_{\gamma}(z)\in J(R)[z]$ and $h_{\gamma}(z)\in R[z]$.
Thus%
\[
(\mathbf{fP})_{\gamma}=\underset{\delta\in\Gamma}{%
{\displaystyle\sum}
}(g_{\delta}(z)+z^{k_{\delta}}h_{\delta}(z))p_{\delta,\gamma}(z)=\underset
{\delta\in\Gamma}{%
{\displaystyle\sum}
}g_{\delta}(z)p_{\delta,\gamma}(z)+h_{\delta}(z)z^{k_{\delta}}p_{\delta
,\gamma}(z)
\]
is in $J(R)[z]+(z^{k_{\gamma}})$, whence $\mathbf{fP}\in\ker(\Phi)$ and
$\mathbf{P}\in\mathcal{M}(X)$ follows. $\square$

\bigskip

\noindent\textbf{4.3.Lemma.}\textit{ If the centre }$Z(R)$\textit{ of the ring
}$R$\textit{ is a field such that }$R/J(R)$\textit{ is finite dimensional over
}$Z(R)$\textit{, then we can exhibit a vector space base of the factor }%
$Z(R)$\textit{-algebra }$\mathcal{N}(X)/\mathcal{I}(X)$\textit{ as}%
\[
\{b_{t}z^{i}\mathbf{E}_{\delta,\gamma}+\mathcal{I}(X)\mid\delta,\gamma
\in\Gamma,k_{\delta,\gamma}\leq i\leq k_{\gamma}-1\text{ and }1\leq
t\leq\lbrack R/J(R):Z(R)]\},
\]
\textit{where the }$b_{t}$\textit{'s are fixed elements of }$R$\textit{ such
that}%
\[
\{b_{t}+J(R)\mid1\leq t\leq\lbrack R/J(R):Z(R)]\}
\]
\textit{is a vector space base of }$R/J(R)$\textit{ over }$Z(R)$\textit{.
Hence}%
\[
\dim_{Z(R)}(\mathcal{N}(X)/\mathcal{I}(X))=[R/J(R):Z(R)]\underset
{\delta,\gamma\in\Gamma}{\sum}(k_{\gamma}-k_{\delta,\gamma}).
\]
\textit{Using }$\Gamma=\{1,2,\ldots,m\}$\textit{ and the assumption that
}$k_{1}\geq k_{2}\geq\ldots\geq k_{m}\geq1$\textit{, we obtain that}%
\[
\dim_{Z(R)}(\mathcal{N}(X)/\mathcal{I}(X))=[R/J(R):Z(R)](k_{1}+3k_{2}%
+\cdots+(2m-1)k_{m}).
\]

\bigskip

\noindent\textbf{Proof.} If $\mathbf{P}=[p_{\delta,\gamma}(z)]$ is in
$\mathcal{N}(X)$, then $p_{\delta,\gamma}(z)\in J(R)[z]+(z^{k_{\delta,\gamma}%
})$ as observed earlier. Thus we have%
\[
p_{\delta,\gamma}(z)=f_{\delta,\gamma}(z)+\left(  \underset{k_{\delta,\gamma
}\leq i\leq k_{\gamma}-1}{\sum}a_{\delta,\gamma,i}z^{i}\right)  +z^{k_{\gamma
}}g_{\delta,\gamma}(z)
\]
for some $f_{\delta,\gamma}(z)\in J(R)[z]$, $a_{\delta,\gamma,i}\in R$ and
$g_{\delta,\gamma}(z)\in R[z]$. In view of the definition of $\mathcal{I}(X)$
we have%
\[
\mathbf{P}+\mathcal{I}(X)=\underset{\delta,\gamma\in\Gamma,k_{\delta,\gamma
}\leq i\leq k_{\gamma}-1}{\sum}(a_{\delta,\gamma,i}z^{i}\mathbf{E}%
_{\delta,\gamma}+\mathcal{I}(X))=
\]%
\[
\underset{\delta,\gamma\in\Gamma,k_{\delta,\gamma}\leq i\leq k_{\gamma}%
-1}{\sum}\underset{t=1}{\overset{[R/J(R):Z(R)]}{\sum}}c_{t}(\delta
,\gamma,i)(b_{t}z^{i}\mathbf{E}_{\delta,\gamma}+\mathcal{I}(X)),
\]
where%
\[
a_{\delta,\gamma,i}+J(R)=\underset{t=1}{\overset{[R/J(R):Z(R)]}{\sum}}%
c_{t}(\delta,\gamma,i)(b_{t}+J(R))
\]
for some $c_{t}(\delta,\gamma,i)\in Z(R)$. Therefore the cosets $b_{t}%
z^{i}\mathbf{E}_{\delta,\gamma}+\mathcal{I}(X)$ generate $\mathcal{N}%
(X)/\mathcal{I}(X)$ over $Z(R)$. It is straightforward to check the
$Z(R)$-linear independence of these cosets. $\square$

\bigskip

\noindent\textbf{4.4.Lemma.}\textit{ The ideal }$zM_{m}(R[z])\vartriangleleft
M_{m}(R[z])$\textit{ is nilpotent modulo }$\mathcal{I}(X)$\textit{, more
precisely we have }$\left(  zM_{m}(R[z])\right)  ^{n}\subseteq\mathcal{I}%
(X)$\textit{, where }$n=\max\{k_{\gamma}\mid\gamma\in\Gamma\}$\textit{. There
is a natural isomorphism between the factor ring }$\mathcal{N}(X)/\left(
\mathcal{N}(X)\cap zM_{m}(R[z])\right)  +\mathcal{I}(X)$\textit{ and the
subring}%
\[
\mathcal{U}(X)=\{U=[\overline{w}_{\delta,\gamma}]\mid\overline{w}%
_{\delta,\gamma}\in R/J(R)\text{ and }\overline{w}_{\delta,\gamma}=0\text{ if
}1\leq k_{\delta}<k_{\gamma}\leq n\}
\]
\textit{of }$M_{m}(R/J(R))$\textit{:}%
\[
\mathcal{N}(X)/\left(  \mathcal{N}(X)\cap zM_{m}(R[z])\right)  +\mathcal{I}%
(X)\cong\mathcal{U}(X)
\]
\textit{and this is an }$(R,R)$\textit{-bimodule isomorphism at the same time.
The ideal}%
\[
\left(  \mathcal{N}(X)\cap zM_{m}(R[z])\right)  +\mathcal{I}(X)/\mathcal{I}%
(X)\vartriangleleft\mathcal{N}(X)/\mathcal{I}(X)
\]
\textit{is nilpotent with }$\left(  \mathcal{N}(X)\cap zM_{m}%
(R[z])+\mathcal{I}(X)/\mathcal{I}(X)\right)  ^{n}=\{0\}$\textit{ and we have
the following isomorphism for the iterated factor:}%
\[
\left(  \mathcal{N}(X)/\mathcal{I}(X)\right)  \!/\!\left(  \mathcal{N}%
(X)\!\cap\!zM_{m}(R[z])\!+\!\mathcal{I}(X)/\mathcal{I}(X)\right)
\!\cong\!\mathcal{N}(X)\!/\!\left(  \mathcal{N}(X)\!\cap\!zM_{m}(R[z])\right)
\!+\!\mathcal{I}(X).
\]

\bigskip

\noindent\textbf{Proof.} Since any entry in the product of the matrices
$\mathbf{Q}_{1},\mathbf{Q}_{2},\ldots,\mathbf{Q}_{n}\in zM_{m}(R[z])$ with
$\mathbf{Q}_{i}=[zq_{\delta,\gamma}^{(i)}(z)]$ is a sum of terms of the form%
\[
zq_{\delta_{1},\gamma_{1}}^{(1)}(z)zq_{\delta_{2},\gamma_{2}}^{(2)}(z)\cdots
zq_{\delta_{n},\gamma_{n}}^{(n)}(z)=z^{n}q_{\delta_{1},\gamma_{1}}%
^{(1)}(z)q_{\delta_{2},\gamma_{2}}^{(2)}(z)\cdots q_{\delta_{n},\gamma_{n}%
}^{(n)}(z),
\]
which is in $(z^{n})$, and since $(z^{n})\subseteq(z^{k_{\gamma}})$ for each
$\gamma\in\Gamma$, we obtain that

\noindent$\left(  zM_{m}(R[z])\right)  ^{n}\subseteq\mathcal{I}(X)$. It
follows that $\left(  \mathcal{N}(X)\cap zM_{m}(R[z])+\mathcal{I}%
(X)/\mathcal{I}(X)\right)  ^{n}=~\{0\}$.

\noindent If $\mathbf{P}=[p_{\delta,\gamma}(z)]$ is in $\mathcal{N}(X)$, then%
\[
p_{\delta,\gamma}(z)=u_{\delta,\gamma}+zf_{\delta,\gamma}(z)+z^{k_{\delta
,\gamma}}(v_{\delta,\gamma}+zg_{\delta,\gamma}(z))
\]
for some $u_{\delta,\gamma}\in J(R)$, $f_{\delta,\gamma}(z)\in J(R)[z]$,
$v_{\delta,\gamma}\in R$ and $g_{\delta,\gamma}(z)\in R[z]$.

\noindent If $1\leq k_{\delta}<k_{\gamma}\leq n$ , then $1\leq k_{\gamma
}-k_{\delta}=k_{\delta,\gamma}$ and%
\[
p_{\delta,\gamma}(z)-u_{\delta,\gamma}=zf_{\delta,\gamma}(z)+z^{k_{\delta
,\gamma}}(v_{\delta,\gamma}+zg_{\delta,\gamma}(z))\in(J(R)[z]+(z^{k_{\delta
,\gamma}}))\cap(zR[z]).
\]
If $1\leq k_{\gamma}\leq k_{\delta}\leq n$, then $k_{\delta,\gamma}=0$,
$z^{k_{\delta,\gamma}}=1$, $(z^{k_{\delta,\gamma}})=R[z]$\ and%
\[
p_{\delta,\gamma}(z)-(u_{\delta,\gamma}+v_{\delta,\gamma})=zf_{\delta,\gamma
}(z)+zg_{\delta,\gamma}(z))\in(J(R)[z]+(z^{k_{\delta,\gamma}}))\cap(zR[z]).
\]
Thus $[w_{\delta,\gamma}]\in M_{m}(R)\cap\mathcal{N}(X)$, $\mathbf{P-}%
[w_{\delta,\gamma}]\in\mathcal{N}(X)\cap zM_{m}(R[z])$ and%
\[
\mathbf{P}+(\left(  \mathcal{N}(X)\cap zM_{m}(R[z])\right)  +\mathcal{I}%
(X))=[w_{\delta,\gamma}]+(\left(  \mathcal{N}(X)\cap zM_{m}(R[z])\right)
+\mathcal{I}(X))
\]
in $\mathcal{N}(X)/\left(  \mathcal{N}(X)\cap zM_{m}(R[z])\right)
+\mathcal{I}(X)$, where%
\[
w_{\delta,\gamma}=\left\{
\begin{array}
[c]{c}%
u_{\delta,\gamma}+v_{\delta,\gamma}\text{ if }1\leq k_{\gamma}\leq k_{\delta
}\leq n\\
u_{\delta,\gamma}\text{ \ \ \ \ \ \ \ \ \ if }1\leq k_{\delta}<k_{\gamma}\leq
n
\end{array}
\right.  .
\]
If $[w_{\delta,\gamma}^{\prime}],[w_{\delta,\gamma}^{\prime\prime}]\in
M_{m}(R)\cap\mathcal{N}(X)$ and%
\[
\lbrack w_{\delta,\gamma}^{\prime}]+\left(  \mathcal{N}(X)\cap zM_{m}%
(R[z])\right)  +\mathcal{I}(X)=[w_{\delta,\gamma}^{\prime\prime}]+\left(
\mathcal{N}(X)\cap zM_{m}(R[z])\right)  +\mathcal{I}(X),
\]
then $w_{\delta,\gamma}^{\prime}+J(R)=w_{\delta,\gamma}^{\prime\prime}+J(R)$
obviously holds in $R/J(R)$\ (for all $\delta,\gamma\in\Gamma$). It follows
that the assignment%
\[
\mathbf{P}+\left(  \mathcal{N}(X)\cap zM_{m}(R[z])\right)  +\mathcal{I}%
(X)\longmapsto\lbrack w_{\delta,\gamma}+J(R)]
\]
is well defined and gives an%
\[
\mathcal{N}(X)/\left(  \mathcal{N}(X)\cap zM_{m}(R[z])\right)  +\mathcal{I}%
(X)\longrightarrow\mathcal{U}(X)
\]
isomorphism. The isomorphism for the iterated factor is obvious. $\square$

\bigskip

\noindent We note that, if $\Gamma=\{1,2,\ldots,m\}$ and $k_{1}>k_{2}%
>\ldots>k_{m}\geq1$, then%
\[
\mathcal{U}(X)=\left[
\begin{array}
[c]{ccccc}%
R/J(R) & R/J(R) & \cdots & R/J(R) & R/J(R)\\
0 & R/J(R) & \ddots & \ddots & R/J(R)\\
\vdots & \ddots & \ddots & \ddots & \vdots\\
0 & \ddots & \ddots & R/J(R) & R/J(R)\\
0 & 0 & \cdots & 0 & R/J(R)
\end{array}
\right]
\]
is an upper triangular\ matrix algebra. In general, if $k_{1}\geq k_{2}%
\geq\ldots\geq k_{m}\geq1$, then $\mathcal{U}(X)$ is a blocked upper
triangular matrix algebra over $R/J(R)$ and the T-ideal of the identities of
$\mathcal{U}(X)$ is determined by Lewin's theorem (see [4]).

\bigskip

\noindent\textbf{4.5.Lemma.}\textit{ For }$\mathbf{P}\in\mathcal{M}%
(X)$\textit{ and }$\mathbf{f}=(f_{\gamma}(z))_{\gamma\in\Gamma}$\textit{ in
}$(R[z])^{\Gamma}$\textit{ the formula}%
\[
\psi_{\mathbf{P}}(\Phi(\mathbf{f}))=\Phi(\mathbf{fP})
\]
\textit{properly defines an }$R$\textit{-endomorphism }$\psi_{\mathbf{P}%
}:M\longrightarrow M$\textit{ of }$_{R}M$\textit{ such that }$\psi
_{\mathbf{P}}\circ\varphi=\varphi\circ\psi_{\mathbf{P}}$\textit{. The
assignment }$\Lambda(\mathbf{P})=\psi_{\mathbf{P}}$\textit{ gives an
}$\mathcal{M}(X)^{\text{op}}\longrightarrow$ Cen$(\varphi)$\textit{
homomorphism of }$Z(R)$\textit{-algebras.}

\bigskip

\noindent\textbf{Proof.} Let $\mathbf{g}\in(R[z])^{\Gamma}$. If $\Phi
(\mathbf{f})=\Phi(\mathbf{g})$, then $\mathbf{f}-\mathbf{g}\in\ker(\Phi)$
implies that $(\mathbf{f}-\mathbf{g})\mathbf{P}\in\ker(\Phi)$, whence
$\Phi(\mathbf{fP})=\Phi(\mathbf{gP})$ follows. Since $\Phi$ is surjective, it
follows that $\psi_{\mathbf{P}}$ is well defined. It is straightforward to
check that%
\[
\psi_{\mathbf{P}}(\Phi(\mathbf{f})+\Phi(\mathbf{g}))=\psi_{\mathbf{P}}%
(\Phi(\mathbf{f}))+\psi_{\mathbf{P}}(\Phi(\mathbf{g}))\text{ and }%
\psi_{\mathbf{P}}(r\Phi(\mathbf{f}))=r\psi_{\mathbf{P}}(\Phi(\mathbf{f}))
\]
for all $\mathbf{f},\mathbf{g}\in(R[z])^{\Gamma}$ and $r\in R$. Thus
$\psi_{\mathbf{P}}$ is an $R$-endomorphism. In view of%
\[
\psi_{\mathbf{P}}(\varphi(\Phi(\mathbf{f}))=\psi_{\mathbf{P}}(z\ast
\Phi(\mathbf{f}))=\psi_{\mathbf{P}}(\Phi(z\mathbf{f}))=\Phi((z\mathbf{f}%
)\mathbf{P})=
\]%
\[
=\Phi(z(\mathbf{fP}))=z\ast\Phi(\mathbf{fP})=z\ast\psi_{\mathbf{P}}%
(\Phi(\mathbf{f}))=\varphi(\psi_{\mathbf{P}}(\Phi(\mathbf{f}))),
\]
the surjectivity of $\Phi$\ gives that $\psi_{\mathbf{P}}\circ\varphi
=\varphi\circ\psi_{\mathbf{P}}$. Clearly,%
\[
\psi_{c\mathbf{P}}=c\psi_{\mathbf{P}}\text{ , }\psi_{\mathbf{P}+\mathbf{Q}%
}=\psi_{\mathbf{P}}+\psi_{\mathbf{Q}}\text{ and }\psi_{\mathbf{PQ}}%
=\psi_{\mathbf{Q}}\circ\psi_{\mathbf{P}}%
\]
ensure that $\Lambda$\ is a homomorphism of $Z(R)$-algebras. We deal only with
the last identity:%
\[
\psi_{\mathbf{PQ}}(\Phi(\mathbf{f}))=\Phi(\mathbf{f}(\mathbf{PQ}%
))=\Phi((\mathbf{fP})\mathbf{Q})=\psi_{\mathbf{Q}}(\Phi(\mathbf{fP}%
))=\psi_{\mathbf{Q}}(\psi_{\mathbf{P}}(\Phi(\mathbf{f}))
\]
proves our claim. $\square$

\bigskip

\noindent\textbf{4.6.Lemma.}\textit{ }$\mathcal{I}(X)\subseteq\ker(\Lambda
)$\textit{ (}$\Lambda$\textit{ is defined in Lemma 4.5). If }$R$\textit{ is a
local ring then }$\mathcal{I}(X)=\ker(\Lambda)$\textit{.}

\bigskip

\noindent\textbf{Proof.} If $\mathbf{P}=[p_{\delta,\gamma}(z)]$ is an element
of $\mathcal{I}(X)$ and $\mathbf{f}=(f_{\gamma}(z))_{\gamma\in\Gamma}$ is an
element of $(R[z])^{\Gamma}$, then $p_{\delta,\gamma}(z)\in
J(R)[z]+(z^{k_{\gamma}})$ implies that%
\[
(\mathbf{fP})_{\gamma}=\underset{\delta\in\Gamma}{%
{\displaystyle\sum}
}f_{\delta}(z)p_{\delta,\gamma}(z)
\]
is in $J(R)[z]+(z^{k_{\gamma}})$\ for all $\gamma\in\Gamma$, whence
$\mathbf{fP}\in\ker(\Phi)$ follows by Proposition 3.1. Since $\psi
_{\mathbf{P}}(\Phi(\mathbf{f}))=\Phi(\mathbf{fP})=0$, we obtain that
$\Lambda(\mathbf{P})=\psi_{\mathbf{P}}=0$, i.e. that $\mathbf{P}\in
\ker(\Lambda)$. Thus the containment is proved.

\noindent If $R$ is a local ring and $\mathbf{P}\in\ker(\Lambda)$, then
$\Lambda(\mathbf{P})=\psi_{\mathbf{P}}=0$ implies that $\psi_{\mathbf{P}}%
(\Phi(\mathbf{f}))=\Phi(\mathbf{fP})=0$ for all $\mathbf{f}\in(R[z])^{\Gamma}%
$. If $\mathbf{1}_{\delta}$ denotes the vector in $(R[z])^{\Gamma}$\ with $1$
in its $\delta$-coordinate and zeros in all other places, then $\mathbf{1}%
_{\delta}\mathbf{P}\in\ker(\Phi)$ and Proposition 3.1 imply that
$p_{\delta,\gamma}(z)\in J(R)[z]+(z^{k_{\gamma}})$. $\square$

\bigskip

\noindent\textbf{4.7.Lemma.}\textit{ If }$\psi\circ\varphi=\varphi\circ\psi
$\textit{ holds for an }$R$\textit{-endomorphism }$\psi:M\longrightarrow
M$\textit{ of }$_{R}M$\textit{, then there exists an }$m\times m$\textit{
matrix }$\mathbf{P}\in\mathcal{M}(X)$\textit{ such that}%
\[
\psi(\Phi(\mathbf{f}))=\Phi(\mathbf{fP})
\]
\textit{for all }$\mathbf{f}=(f_{\gamma}(z))_{\gamma\in\Gamma}$\textit{ in
}$(R[z])^{\Gamma}$\textit{.}

\bigskip

\noindent\textbf{Proof.} Since $\Phi:(R[z])^{\Gamma}\longrightarrow M$ is
surjective, for each $\delta\in\Gamma$\ we can find an element $\mathbf{p}%
_{\delta}=(p_{\delta,\gamma}(z))_{\gamma\in\Gamma}$ in $(R[z])^{\Gamma}$ such
that $\Phi(\mathbf{p}_{\delta})=\psi(x_{\delta,1})$. For the $m\times m$
matrix $\mathbf{P}=[p_{\delta,\gamma}(z)]$ we have%
\[
\psi(\Phi(\mathbf{f}))=\underset{\delta\in\Gamma}{%
{\displaystyle\sum}
}\psi(f_{\delta}(z)\ast x_{\delta,1})=\underset{\delta\in\Gamma}{%
{\displaystyle\sum}
}f_{\delta}(z)\ast\psi(x_{\delta,1})=\underset{\delta\in\Gamma}{%
{\displaystyle\sum}
}f_{\delta}(z)\ast\Phi(\mathbf{p}_{\delta})=
\]%
\[
=\underset{\delta\in\Gamma}{%
{\displaystyle\sum}
}\Phi(f_{\delta}(z)\mathbf{p}_{\delta})=\Phi(\underset{\delta\in\Gamma}{%
{\displaystyle\sum}
}f_{\delta}(z)\mathbf{p}_{\delta})=\Phi(\mathbf{fP})
\]
for all $\mathbf{f}\in(R[z])^{\Gamma}$. Since $\mathbf{f}\in\ker(\Phi)$
implies that $\Phi(\mathbf{fP})=\psi(\Phi(\mathbf{f}))=0$, we obtain that
$\mathbf{P}\in\mathcal{M}(X)$. $\square$

\bigskip

\noindent\textbf{4.8.Theorem.}\textit{ Let }$\varphi:M\longrightarrow
M$\textit{ be a nilpotent }$R$\textit{-endomorphism of the finitely generated
semisimple left }$R$\textit{-module }$_{R}M$\textit{. Then }$\Lambda
:\mathcal{M}(X)^{\text{op}}\longrightarrow$ Cen$(\varphi)$\textit{ (defined in
Lemma 4.5) is a surjective homomorphism of }$Z(R)$\textit{-algebras, where the
centralizer }Cen$(\varphi)$\textit{\ is a }$Z(R)$\textit{-subalgebra of
}Hom$_{R}(M,M)$\textit{ and }$m=\dim_{R}(\ker(\varphi))$\textit{.}

\bigskip

\noindent\textbf{Proof.} Lemma 4.5 ensures that $\Lambda:\mathcal{M}%
(X)^{\text{op}}\longrightarrow$ Cen$(\varphi)$ is a homomorphism of
$Z(R)$-algebras. The surjectivity of $\Lambda$\ follows from Lemma 4.7. To
conclude the proof it suffices to note that $m=\left\vert \Gamma\right\vert
=\dim_{R}(\ker(\varphi))$. $\square$

\bigskip

\noindent\textbf{4.9.Corollary.}\textit{ Let }$\varphi:M\longrightarrow
M$\textit{ be a nilpotent }$R$\textit{-endomorphism of the finitely generated
semisimple left }$R$\textit{-module }$_{R}M$\textit{. Then }Cen$(\varphi
)$\textit{ satisfies all of the polynomial identities (with coefficients in
}$Z(R)$\textit{) of }$M_{m}^{\text{op}}(R[z])$\textit{. If }$R$\textit{\ is
commutative, then }Cen$(\varphi)$\textit{ satisfies the standard identity
}$S_{2m}=0$\textit{ of degree }$2m$\textit{ by the Amitsur-Levitzki theorem.}

\bigskip

\noindent\textbf{4.10.Theorem.}\textit{ Let }$R$\textit{\ be a local ring and
}$\varphi:M\longrightarrow M$\textit{ a nilpotent }$R$\textit{-endomorphism of
the finitely generated semisimple left }$R$\textit{-module }$_{R}M$\textit{.
Then the centralizer }Cen$(\varphi)$\textit{\ of }$\varphi$\textit{\ is
isomorphic to the opposite of the factor }$\mathcal{N}(X)/\mathcal{I}%
(X)$\textit{ as }$Z(R)$\textit{-algebras:}%
\[
\text{Cen}(\varphi)\cong\left(  \mathcal{N}(X)/\mathcal{I}(X)\right)
^{\text{op}}\cong\mathcal{N}^{\text{op}}(X)/\mathcal{I}(X).
\]
\textit{If }$f_{i}=0$\textit{ are polynomial identities of the }%
$Z(R)$\textit{-subalgebra }$\mathcal{U}^{\text{op}}(X)$\textit{ of }%
$M_{m}^{\text{op}}(R/J(R))$\textit{ with }$f_{i}\in K\langle x_{1}%
,\ldots,x_{r}\rangle$\textit{, }$1\leq i\leq n$\textit{, then }$f_{1}%
f_{2}\cdots f_{n}=0$\textit{ is an identity of }Cen$(\varphi)$\textit{ (here
}$\varphi^{n}=0\neq\varphi^{n-1}$\textit{ and }$m=\dim_{R}(\ker(\varphi
))$\textit{).}

\bigskip

\noindent\textbf{Proof.} Theorem 4.8 ensures that Cen$(\varphi)\cong
\mathcal{M}(X)^{\text{op}}/\ker(\Lambda)$ as $Z(R)$-algebras. In order to
prove the desired isomorphism, it suffices to note that for a local ring $R$
we have $\mathcal{M}(X)=\mathcal{N}(X)$ and $\ker(\Lambda)=\mathcal{I}(X)$ by
Lemmas 4.2 and 4.6 respectively. Now Lemma 4.4 ensures that%
\[
L=\left(  \mathcal{N}(X)\cap zM_{m}(R[z])\right)  +\mathcal{I}(X)/\mathcal{I}%
(X)\vartriangleleft\mathcal{N}(X)/\mathcal{I}(X)
\]
can be viewed as an ideal of Cen$(\varphi)$ such that $L^{n}=\{0\}$ and%
\[
\text{Cen}(\varphi)/L\cong(\mathcal{N}^{\text{op}}(X)/\mathcal{I}(X))\diagup
L\cong\mathcal{N}^{\text{op}}(X)\diagup\left(  \mathcal{N}(X)\cap
zM_{m}(R[z])\right)  +\mathcal{I}(X)\cong\mathcal{U}^{\text{op}}(X).
\]
It follows that $f_{i}=0$ is an identity of Cen$(\varphi)/L$. Thus
$f_{i}(v_{1},\ldots,v_{r})\in L$\ for all $v_{1},\ldots,v_{r}\in$%
Cen$(\varphi)$, whence we obtain that $f_{1}f_{2}\cdots f_{n}=0$ is an
identity of Cen$(\varphi)$. $\square$

\bigskip

\noindent\textbf{4.11.Theorem.}\textit{ Let }$R$\textit{\ be a local ring such
that }$Z(R)$\textit{ is a field and }$R/J(R)$\textit{ is finite dimensional
over }$Z(R)$\textit{. If }$\varphi:M\longrightarrow M$\textit{ is a nilpotent
}$R$\textit{-endomorphism of the finitely generated semisimple left }%
$R$\textit{-module }$_{R}M$\textit{, then}%
\[
\dim_{Z(R)}(\text{Cen}(\varphi))=[R/J(R):Z(R)](k_{1}+3k_{2}+\cdots
+(2m-1)k_{m}),
\]
\textit{where }$k_{1}\geq k_{2}\geq\ldots\geq k_{m}\geq1$\textit{ are the
sizes of the blocks of the nilpotent Jordan normal base }$X$\textit{ with
respect to }$\varphi$\textit{.}

\bigskip

\noindent\textbf{Proof.} By Theorem 4.10 we have%
\[
\text{Cen}(\varphi)\cong\left(  \mathcal{N}(X)/\mathcal{I}(X)\right)
^{\text{op}},
\]
and since%
\[
\dim_{Z(R)}\left(  \mathcal{N}(X)/\mathcal{I}(X)\right)  ^{\text{op}}%
=\dim_{Z(R)}(\mathcal{N}(X)/\mathcal{I}(X)),
\]
the result follows from Lemma 4.3. $\square$

\bigskip

\noindent In the nilpotent case, Theorem 4.11 generalizes the formula for the
dimension of the centralizer Cen$(A)$ of a matrix $A\in M_{n}(K)$ over a field
$K$ (see [9,10]).

\bigskip

\noindent5. FURTHER\ PROPERTIES\ OF\ THE\ CENTRALIZERS

\bigskip

\noindent\textbf{5.1.Theorem.}\textit{ If }$_{R}M$\textit{ is semisimple and
}$\varphi:M\longrightarrow M$\textit{ is an indecomposable nilpotent element
of the ring }Hom$_{R}(M,M)$\textit{, then the following are equivalent.}

\begin{enumerate}
\item $\psi\in$ Cen$(\varphi)$\textit{.}

\item \textit{We can find an }$R$\textit{-generating set }$\{y_{j}\in
M\mid1\leq j\leq d\}$\textit{ of }$_{R}M$\textit{ and elements }$a_{1}%
,a_{2},\ldots,a_{n}$\textit{ in }$R$\textit{ such that}%
\[
a_{1}y_{j}+a_{2}\varphi(y_{j})+\cdots+a_{n}\varphi^{n-1}(y_{j})=\psi(y_{j})
\]
\textit{and}%
\[
a_{1}\varphi(y_{j})+a_{2}\varphi(\varphi(y_{j}))+\cdots+a_{n}\varphi
^{n-1}(\varphi(y_{j}))=\psi(\varphi(y_{j}))
\]
\textit{for all }$1\leq j\leq d$\textit{.}
\end{enumerate}

\bigskip

\noindent\textbf{Proof.}

\noindent(1)$\Longrightarrow$(2): Obviously, if $\psi\in$ Cen$(\varphi)$ then
the first identity implies the second one. Proposition 2.3 ensures the
existence of a nilpotent Jordan normal base $\{x_{i}\mid1\leq i\leq n\}$\ of
$_{R}M$\ with respect to $\varphi$ consisting of one block. Clearly,
$\underset{1\leq i\leq n}{\oplus}Rx_{i}=M$ implies that%
\[
\psi(x_{1})=a_{1}x_{1}+a_{2}x_{2}+\cdots+a_{n}x_{n}=a_{1}x_{1}+a_{2}%
\varphi(x_{1})+\cdots+a_{n}\varphi^{n-1}(x_{1})
\]
for some $a_{1},a_{2},\ldots,a_{n}\in R$. Thus%
\[
\psi(x_{i})=\psi(\varphi^{i-1}(x_{1}))=\varphi^{i-1}(\psi(x_{1}))=\varphi
^{i-1}(a_{1}x_{1}+a_{2}\varphi(x_{1})+\cdots+a_{n}\varphi^{n-1}(x_{1}))=
\]%
\[
=\!a_{1}\varphi^{i-1}(x_{1})\!+\!a_{2}\varphi(\varphi^{i-1}(x_{1}%
))\!+\cdots+\!a_{n}\varphi^{n-1}(\varphi^{i-1}(x_{1}))\!=\!a_{1}%
x_{i}\!+\!a_{2}\varphi(x_{i})\!+\cdots+\!a_{n}\varphi^{n-1}(x_{i})
\]
for all $1\leq i\leq n$.

\noindent(2)$\Longrightarrow$(1): Since we have%
\[
\varphi(\psi(y_{j}))=\varphi(a_{1}y_{j}+a_{2}\varphi(y_{j})+\cdots
+a_{n}\varphi^{n-1}(y_{j}))=
\]%
\[
=a_{1}\varphi(y_{j})+a_{2}\varphi(\varphi(y_{j}))+\cdots+a_{n}\varphi
^{n-1}(\varphi(y_{j}))=\psi(\varphi(y_{j}))
\]
for all $1\leq j\leq d$, the implication is proved. $\square$

\bigskip

\noindent\textbf{5.2.Corollary.}\textit{ If }$R$\textit{ is commutative,
}$_{R}M$\textit{ is semisimple and }$\varphi:M\longrightarrow M$\textit{ is an
indecomposable nilpotent element of the ring }Hom$_{R}(M,M)$\textit{, then the
following are equivalent.}

\begin{enumerate}
\item $\psi\in$ Cen$(\varphi)$\textit{.}

\item \textit{We can find elements }$a_{1},a_{2},\ldots,a_{n}$\textit{ in }%
$R$\textit{ such that}%
\[
a_{1}u+a_{2}\varphi(u)+\cdots+a_{n}\varphi^{n-1}(u)=\psi(u)
\]
\textit{for all }$u\in M$\textit{. In other words, }$\psi$\textit{\ is a
polynomial of }$\varphi$\textit{.}
\end{enumerate}

\bigskip

\noindent\textbf{Proof.} It suffices to prove that if $\underset{1\leq j\leq
d}{\sum}Ry_{j}=M$ and%
\[
a_{1}y_{j}+a_{2}\varphi(y_{j})+\cdots+a_{n}\varphi^{n-1}(y_{j})=\psi(y_{j})
\]
holds for all $1\leq j\leq d$, then we have%
\[
a_{1}u+a_{2}\varphi(u)+\cdots+a_{n}\varphi^{n-1}(u)=\psi(u)
\]
for all $u\in M$. Since $u=b_{1}y_{1}+b_{2}y_{2}+\cdots+b_{d}y_{d}$ for some
$b_{1},b_{2},\ldots,b_{d}\in R$ and $b_{j}a_{i}=a_{i}b_{j}$, we obtain that%
\[
\psi(u)=\underset{1\leq j\leq d}{\sum}b_{j}\psi(y_{j})=\underset{1\leq j\leq
d}{\sum}b_{j}(a_{1}y_{j}+a_{2}\varphi(y_{j})+\cdots+a_{n}\varphi^{n-1}%
(y_{j}))=
\]%
\[
=a_{1}\left(  \underset{1\leq j\leq d}{\sum}b_{j}y_{j}\right)  +a_{2}%
\varphi\left(  \underset{1\leq j\leq d}{\sum}b_{j}y_{j}\right)  +\cdots
+a_{n}\varphi^{n-1}\left(  \underset{1\leq j\leq d}{\sum}b_{j}y_{j}\right)  =
\]%
\[
=a_{1}u+a_{2}\varphi(u)+\cdots+a_{n}\varphi^{n-1}(u).\ \square
\]

\bigskip

\noindent\textbf{5.3.Theorem.}\textit{ Let }$R$\textit{\ be a local ring. If
}$\varphi:M\longrightarrow M$\textit{ is a nilpotent }$R$\textit{-endomorphism
of the finitely generated semisimple left }$R$\textit{-module }$_{R}M$\textit{
and}

\noindent$\sigma\in$ Hom$_{R}(M,M)$\textit{ is arbitrary, then the following
are equivalent.}

\begin{enumerate}
\item Cen$(\varphi)\subseteq$ Cen$(\sigma)$\textit{.}

\item \textit{We can find an }$R$\textit{-generating set }$\{y_{j}\in
M\mid1\leq j\leq d\}$\textit{ of }$_{R}M$\textit{ and elements }$a_{1}%
,a_{2},\ldots,a_{n}$\textit{ in }$R$\textit{ such that}%
\[
a_{1}\psi(y_{j})+a_{2}\varphi(\psi(y_{j}))+\cdots+a_{n}\varphi^{n-1}%
(\psi(y_{j}))=\sigma(\psi(y_{j}))
\]
\textit{for all }$1\leq j\leq d$\textit{ and all }$\psi\in$ Cen$(\varphi
)$\textit{.}
\end{enumerate}

\bigskip

\noindent\textbf{Proof.}

\noindent(1)$\Longrightarrow$(2): Obviously, if Cen$(\varphi)\subseteq$
Cen$(\sigma)$ then%
\[
a_{1}y_{j}+a_{2}\varphi(y_{j})+\cdots+a_{n}\varphi^{n-1}(y_{j})=\sigma(y_{j})
\]
implies that%
\[
a_{1}\psi(y_{j})+a_{2}\varphi(\psi(y_{j}))+\cdots+a_{n}\varphi^{n-1}%
(\psi(y_{j}))=\sigma(\psi(y_{j}))
\]
for all $\psi\in$ Cen$(\varphi)$. Theorem 2.1 ensures the existence of a
nilpotent Jordan normal base $\{x_{\gamma,i}\mid\gamma\in\Gamma,1\leq i\leq
k_{\gamma}\}$ of $_{R}M$\ with respect to $\varphi$. Consider the natural
projection $\varepsilon_{\delta}:M\longrightarrow N_{\delta}^{\prime}$
corresponding to the direct sum $M=N_{\delta}^{\prime}\oplus N_{\delta
}^{\prime\prime}$ (see the proof of 2.3), where
\[
N_{\delta}^{\prime}=\underset{1\leq i\leq k_{\delta}}{\oplus}Rx_{\delta
,i}\text{ and }N_{\delta}^{\prime\prime}=\underset{\gamma\in\Gamma
\setminus\{\delta\},1\leq i\leq k_{\gamma}}{\oplus}Rx_{\gamma,i}\text{ }.
\]
Then $\varepsilon_{\delta}\in$ Cen$(\varphi)$, whence $\varepsilon_{\delta}%
\in$ Cen$(\sigma)$ follows for all $\delta\in\Gamma$. Thus im$(\varepsilon
_{\delta})=N_{\delta}^{\prime}$ and%
\[
\sigma:\text{im}(\varepsilon_{\delta})\longrightarrow\text{im}(\varepsilon
_{\delta})
\]
implies that%
\[
\sigma(x_{\delta,1})=\underset{1\leq i\leq k_{\delta}}{%
{\displaystyle\sum}
}a_{\delta,i}x_{\delta,i}=h_{\delta}(z)\ast x_{\delta,1}%
\]
for some $h_{\delta}(z)=a_{\delta,1}+a_{\delta,2}z+\cdots+a_{\delta,k_{\delta
}}z^{k_{\delta}-1}$ in $R[z]$. Since $\varphi\in$ Cen$(\sigma)$ implies that
$\sigma\in$ Cen$(\varphi)$, it follows that%
\[
\sigma(\Phi(\mathbf{f}))=\underset{\gamma\in\Gamma}{%
{\displaystyle\sum}
}\sigma(f_{\gamma}(z)\ast x_{\gamma,1})=\underset{\gamma\in\Gamma}{%
{\displaystyle\sum}
}f_{\gamma}(z)\ast\sigma(x_{\gamma,1})=\underset{\gamma\in\Gamma}{%
{\displaystyle\sum}
}f_{\gamma}(z)\ast(h_{\gamma}(z)\ast x_{\gamma,1})=
\]%
\[
=\underset{\gamma\in\Gamma}{%
{\displaystyle\sum}
}(f_{\gamma}(z)h_{\gamma}(z))\ast x_{\gamma,1}=\Phi(\mathbf{fH}),
\]
where $\mathbf{f}\in(R[z])^{\Gamma}$ and $\mathbf{H}=\underset{\gamma\in
\Gamma}{\sum}h_{\gamma}(z)\mathbf{E}_{\gamma,\gamma}$ is a $\Gamma\times
\Gamma$ diagonal matrix in $\mathcal{M}(X)$ (note that $\mathbf{H\in
}\mathcal{M}(X)$ is a consequence of $\sigma(\Phi(\mathbf{f}))=\Phi
(\mathbf{fH})$). In view of Theorem 4.8, the containment Cen$(\varphi
)\subseteq$ Cen$(\sigma)$ is equivalent to the condition that $\sigma\circ
\psi_{\mathbf{P}}=\psi_{\mathbf{P}}\circ\sigma$ for all $\mathbf{P}%
\in\mathcal{M}(X)$. Consequently, we obtain that Cen$(\varphi)\subseteq
$\ Cen$(\sigma)$ is equivalent to the following:%
\[
\Phi(\mathbf{fPH})=\sigma(\Phi(\mathbf{fP})))=\sigma(\psi_{\mathbf{P}}%
(\Phi(\mathbf{f})))=\psi_{\mathbf{P}}(\sigma(\Phi(\mathbf{f})))=\psi
_{\mathbf{P}}(\Phi(\mathbf{fH}))=\Phi(\mathbf{fHP})
\]
for all $\mathbf{f}\in(R[z])^{\Gamma}$ and $\mathbf{P}\in\mathcal{M}(X)$. Thus
Cen$(\varphi)\subseteq$ Cen$(\sigma)$ implies that $\Phi(\mathbf{f}%
(\mathbf{PH}-\mathbf{HP}))=0$, i.e. that $\mathbf{f}(\mathbf{PH}%
-\mathbf{HP})\in\ker(\Phi)$ for all $\mathbf{f}$ and $\mathbf{P}$.

\noindent Now we use that $R$\ is a local ring. If $\alpha\in\Gamma$ is an
index such that%
\[
k_{\alpha}=\max\{k_{\gamma}\mid\gamma\in\Gamma\}=n,
\]
then $\mathbf{E}_{\alpha,\delta}\in\mathcal{M}(X)$ for all $\delta\in\Gamma$
(see Remark 4.1). Take $\mathbf{e}=(1)_{\gamma\in\Gamma}$ and $\mathbf{P}%
=\mathbf{E}_{\alpha,\delta}$, then the $\delta$-coordinate of%
\[
\mathbf{e}(\mathbf{E}_{\alpha,\delta}\mathbf{H}-\mathbf{HE}_{\alpha,\delta
})=(h_{\delta}(z)-h_{\alpha}(z))\mathbf{eE}_{\alpha,\delta}%
\]
is $h_{\delta}(z)-h_{\alpha}(z)$. Since%
\[
(h_{\delta}(z)-h_{\alpha}(z))\mathbf{eE}_{\alpha,\delta}\in\ker(\Phi
)=\underset{\gamma\in\Gamma}{%
{\displaystyle\coprod}
}J(R)[z]+(t^{k_{\gamma}}),
\]
we obtain that $h_{\delta}(z)-h_{\alpha}(z)\in J(R)[z]+(z^{k_{\delta}})$. Thus%
\[
\sigma(x_{\delta,1})=h_{\delta}(z)\ast x_{\delta,1}=h_{\alpha}(z)\ast
x_{\delta,1}%
\]
for all $\delta\in\Gamma$. It follows that%
\[
\sigma(x_{\gamma,i})=\sigma(\varphi^{i-1}(x_{\gamma,1}))=\varphi^{i-1}%
(\sigma(x_{\gamma,1}))=\varphi^{i-1}(h_{\alpha}(z)\ast x_{\gamma,1})=
\]%
\[
=h_{\alpha}(z)\ast\varphi^{i-1}(x_{\gamma,1})=h_{\alpha}(z)\ast x_{\gamma
,i}=a_{1}x_{\gamma,i}+a_{2}\varphi(x_{\gamma,i})+\cdots+a_{n}\varphi
^{n-1}(x_{\gamma,i}),
\]
where $h_{\alpha}(z)=a_{1}+a_{2}z+\cdots+a_{n}z^{n-1}$.

\noindent(2)$\Longrightarrow$(1): ($R$ is an arbitrary ring) Since $1_{M}\in$
Cen$(\varphi)$, we obtain that%
\[
a_{1}y_{j}+a_{2}\varphi(y_{j})+\cdots+a_{n}\varphi^{n-1}(y_{j})=\sigma(y_{j})
\]
for all $1\leq j\leq d$. If $\psi\in$ Cen$(\varphi)$, then%
\[
\psi(\sigma(y_{j}))=\psi(a_{1}y_{j}+a_{2}\varphi(y_{j})+\cdots+a_{n}%
\varphi^{n-1}(y_{j}))=
\]%
\[
=a_{1}\psi(y_{j})+a_{2}\varphi(\psi(y_{j}))+\cdots+a_{n}\varphi^{n-1}%
(\psi(y_{j}))=\sigma(\psi(y_{j}))
\]
for all $1\leq j\leq d$, whence $\psi\circ\sigma=\sigma\circ\psi$ follows.
Thus Cen$(\varphi)\subseteq$ Cen$(\sigma)$. $\square$

\bigskip

\noindent6. THE\ CENTRALIZER\ OF\ AN ARBITRARY LINEAR MAP

\bigskip

\noindent If the field $K$ is arbitrary and $E$ is an extension of $K$, we may
assume that $M_{n}(K)\subseteq M_{n}(E)$. Let us denote the centralizers of
$A\in M_{n}(K)$ in $M_{n}(K)$ and in $M_{n}(E)$ by Cen$_{K}(A)$ and
Cen$_{E}(A)$, respectively. Since
\[
M_{n}(E)\cong E\otimes_{K}M_{n}(K),
\]
we obtain that
\[
\text{Cen}_{E}(A)\cong E\otimes_{K}\text{Cen}_{K}(A).
\]
If the field $K$ is infinite, then the $K$-algebra $S$ and the $E$-algebra
$E\otimes_{K}S$ have the same polynomial identities. If $T_{K}(S)\subseteq
K\langle x_{1},\ldots,x_{n},...\rangle$ and
\[
T_{E}(E\otimes_{K}S)\subseteq E\otimes_{K}K\langle x_{1},\ldots,x_{n}%
,...\rangle\cong E\langle x_{1},\ldots,x_{n},...\rangle
\]
are the T-ideals of the $K$-algebra $S$ and the $E$-algebra $E\otimes_{K}S$,
respectively, then
\[
T(E\otimes_{K}S)=E\otimes_{K}T(S).
\]
If the field $K$ is finite, this holds for the multilinear identities only.
Hence the information on the (at least multilinear) polynomial identities of
Cen$(A)$ for $K$ arbitrary can be derived from the case when $K$ is
algebraically closed.

\noindent If $\{\lambda_{1},\lambda_{2},\ldots,\lambda_{r}\}$ is the set of
all eigenvalues of $A$, then Cen$(A)$\ is isomorphic to the direct product of
the centralizers Cen$(A_{i})$, where $A_{i}$ denotes the block diagonal matrix
consisting of all Jordan blocks of $A$ having eigenvalue $\lambda_{i}$ in the
diagonal. The number of the diagonal blocks in $A_{i}$ is $\dim(\ker
(A_{i}-\lambda_{i}I_{i}))$, and the size of $A_{i}$ is $d_{i}\times d_{i}$,
where $d_{i}$ is the multiplicity of the root $\lambda_{i}$\ in the
characteristic polynomial of $A$. Since%
\[
\text{Cen}(A_{i})=\text{Cen}(A_{i}-\lambda_{i}I_{i})
\]
and $A_{i}-\lambda_{i}I_{i}$ is nilpotent in $M_{d_{i}}(K)$, we shall consider
the case of a nilpotent matrix. We note that a multiplicative ($K$ vector
space) base of Cen$(A)$ was constructed in [8], where it was also proved that
Cen$(A)$ is uniquely determined by the block structure of the Jordan normal
form of $A$.

\bigskip

\noindent\textbf{6.1.Theorem.}\textit{ Let }$A\in M_{d}(K)$\textit{ be a
nilpotent matrix and let }$J($Cen$(A))$\textit{ be the Jacobson radical of
}Cen$(A)$\textit{. The elementary Jordan matrices in the canonical Jordan form
of }$A$\textit{ are indexed by the elements of }$\Gamma=\{1,2,\ldots
,m\}$\textit{ and we have }$k_{1}\geq k_{2}\geq\ldots\geq k_{m}\geq1$\textit{
for the sizes of these elementary Jordan blocks. Then }%
\[
\text{Cen}(A)/J(\text{Cen}(A))\cong M_{p_{1}}(K)\oplus\cdots\oplus M_{p_{u}%
}(K)
\]
\textit{for some }$u$\textit{, where }$p_{e}$\textit{ is the number of
elementary Jordan matrices of size }$e\times e$\textit{ and }$M_{p_{e}%
}(K)=\{0\}$\textit{ if }$p_{e}=0$\textit{. The index of nilpotency of }%
$J($Cen$(A))$\textit{ is bounded from above by }$nv$\textit{, where }%
$n=\max\{k_{i}\mid i\in\Gamma\}$\textit{ and }$v$\textit{ is the number of
different sizes.}

\bigskip

\noindent\textbf{Proof.} The matrix $A$ can be considered as a nilpotent
$K$-linear map of the vector space $K^{d}$. The Jordan normal form of $A$
provides a nilpotent Jordan normal base in $K^{d}$\ with block sizes
$k_{1}\geq k_{2}\geq\ldots\geq k_{m}\geq1$. Now $k_{i,j}=k_{j}-k_{i}$ when
$1\leq k_{i}<k_{j}\leq n$ and $k_{i,j}=0$ otherwise. The application of
Theorem 4.10 gives an isomorphism Cen$(A)\cong\left(  \mathcal{N}%
(X)/\mathcal{I}(X)\right)  ^{\text{op}}\cong\mathcal{N}^{\text{op}%
}(X)/\mathcal{I}(X)$ of $K$-algebras, where%
\[
\mathcal{N}(X)=\{\mathbf{P}\in M_{m}(K[z])\mid\mathbf{P}=[p_{i,j}(z)]\text{
and }z^{k_{i}}p_{i,j}(z)\in(z^{k_{j}})\text{ for all }1\leq i,j\leq m\},
\]%
\[
\mathcal{I}(X)=\{\mathbf{P}\in M_{m}(K[z])\mid\mathbf{P}=[p_{i,j}(z)]\text{
and }p_{i,j}(z)\in(z^{k_{j}})\text{ for all }1\leq i,j\leq m\}
\]
and $z^{k_{i}}p_{i,j}(z)\in(z^{k_{j}})$ is equivalent to $p_{i,j}%
(z)\in(z^{k_{i,j}})$ (see the three sentences preceding Remark 4.1). Let
$T_{i}=K[z]/(z^{k_{i}})$, and denote by the same symbol $z$ the element
$z+(z^{k_{i}})$ of $T_{i}$. Take $l_{i,j}=k_{j,i}$ for all $i,j\in\Gamma$ and
consider the set%
\[
\mathcal{C}_{A}=\left[
\begin{matrix}
z^{l_{1,1}}T_{1} & z^{l_{1,2}}T_{1} & \cdots & z^{l_{1,m}}T_{1}\\
z^{l_{2,1}}T_{2} & z^{l_{2,2}}T_{2} & \cdots & z^{l_{2,m}}T_{2}\\
\vdots & \vdots & \ddots & \vdots\\
z^{l_{m,1}}T_{m} & z^{l_{m,2}}T_{m} & \cdots & z^{l_{m,m}}T_{m}%
\end{matrix}
\right]  =\sum_{i,j=1}^{m}z^{l_{i,j}}T_{i}E_{i,j}%
\]
of $m\times m$ matrices, where the $E_{i,j}$'s are the usual matrix units in
$M_{m}(K)$. It is straightforward to see that the natural matrix addition and
multiplication give a $K$-algebra structure on $\mathcal{C}_{A}$. Using a
matrix $\mathbf{P}=[p_{i,j}(z)]$ in $\mathcal{N}(X)$, the map%
\[
\mathbf{P}+\mathcal{I}(X)\longmapsto\lbrack p_{i,j}(z)+(z^{k_{j}})]^{\top}%
\]
(here $^{\top}$ denotes the transpose) is well defined and provides an
$\mathcal{N}^{\text{op}}(X)/\mathcal{I}(X)\rightarrow\mathcal{C}_{A}$
isomorphism of $K$-algebras, whence%
\[
\text{Cen}(A)\cong\mathcal{C}_{A}%
\]
can be derived. Recall that the Jacobson radical of a finite dimensional
algebra is equal to the maximal nilpotent ideal of the algebra. Since
$k_{1}\geq k_{2}\geq\ldots\geq k_{m}$, the $K[z]$-module%
\[
\mathcal{T}_{A}=\left[
\begin{matrix}
T_{1} & T_{1} & \cdots & T_{1}\\
T_{2} & T_{2} & \cdots & T_{2}\\
\vdots & \vdots & \ddots & \vdots\\
T_{m} & T_{m} & \cdots & T_{m}%
\end{matrix}
\right]  =\sum_{i,j=1}^{m}T_{i}E_{i,j}%
\]
satisfies $z^{k_{1}}\!\mathcal{T}_{A}\!\!=\!\!\{0\}$. The intersection
$I\!\!=\!\!z\mathcal{T}_{A}\!\cap\mathcal{C}_{A}$ is an ideal of
$\mathcal{C}_{A}$ and $I^{n}\!\!=\!\!I^{k_{1}}\!\!=\!\!\{0\}$. Hence
$I\subseteq J(\mathcal{C}_{A})$. Since $l_{i,j}=k_{i}-k_{j}$ when $i<j$ and
$l_{i,j}=0$ if $i>j$ or $k_{i}=k_{j}$, we obtain that%
\[
\mathcal{C}_{A}/I=\left[
\begin{matrix}
K & \cdots & K & 0 & \cdots & 0 &  & 0 & \cdots & 0\\
\vdots & \ddots & \vdots & \vdots & \ddots & \vdots & \cdots & \vdots & \ddots
& \vdots\\
K & \cdots & K & 0 & \cdots & 0 &  & 0 & \cdots & 0\\
K & \cdots & K & K & \cdots & K &  & 0 & \cdots & 0\\
\vdots & \ddots & \vdots & \vdots & \ddots & \vdots & \cdots & \vdots & \ddots
& \vdots\\
K & \cdots & K & K & \cdots & K &  & 0 & \cdots & 0\\
& \vdots &  &  & \vdots &  & \ddots &  & \vdots & \\
K & \cdots & K & K & \cdots & K &  & K & \cdots & K\\
\vdots & \ddots & \vdots & \vdots & \ddots & \vdots & \cdots & \vdots & \ddots
& \vdots\\
K & \cdots & K & K & \cdots & K &  & K & \cdots & K
\end{matrix}
\right]  .
\]
Here the diagonal blocks
\[
\left[
\begin{matrix}
K & \cdots & K &  &  &  &  &  &  & \\
\vdots & \ddots & \vdots &  &  &  &  &  &  & \\
K & \cdots & K &  &  &  &  &  &  & \\
&  &  & K & \cdots & K &  &  &  & \\
&  &  & \vdots & \ddots & \vdots &  &  &  & \\
&  &  & K & \cdots & K &  &  &  & \\
&  &  &  &  &  & \ddots &  &  & \\
&  &  &  &  &  &  & K & \cdots & K\\
&  &  &  &  &  &  & \vdots & \ddots & \vdots\\
&  &  &  &  &  &  & K & \cdots & K
\end{matrix}
\right]
\]
are matrix algebras of size $p_{1}\times p_{1},p_{2}\times p_{2},\ldots
,p_{u}\times p_{u}$. The number of blocks is equal to the number $v$ of
different sizes of elementary Jordan matrices in the canonical Jordan form of
$A$. Hence the lower triangular part of $\mathcal{C}_{A}/I$%
\[
\left[
\begin{matrix}
\phantom{0} &  &  &  &  &  &  &  &  & \\
\phantom{\vdots} &  &  &  &  &  &  &  &  & \\
\phantom{0} &  &  &  &  &  &  &  &  & \\
K & \cdots & K &  &  &  &  &  &  & \\
\vdots & \ddots & \vdots &  &  &  &  &  &  & \\
K & \cdots & K &  & \ddots &  &  &  &  & \\
& \vdots &  & \ddots &  &  &  &  &  & \\
K & \cdots & K &  & K & \cdots & K &  &  & \\
\vdots & \ddots & \vdots & \cdots & \vdots & \ddots & \vdots &  &  & \\
K & \cdots & K &  & K & \cdots & K & \phantom{0} & \phantom{\cdots} &
\phantom{0}
\end{matrix}
\right]
\]
consists of $v\times v$ lower triangular block matrices, is nilpotent of index
$v$ and is equal to the radical of $\mathcal{C}_{A}/I$. Hence $(J(\mathcal{C}%
_{A})^{v})^{n}\subseteq I^{n}=\{0\}$ and the class of nilpotency of
$J(\mathcal{C}_{A})$ is bounded by $nv$. Clearly
\[
\mathcal{C}_{A}/J(\mathcal{C}_{A})\cong M_{p_{1}}(K)\oplus\cdots\oplus
M_{p_{u}}(K).\ \square
\]

\bigskip

\noindent We note that when $K$ is of characteristic $0$, algebras of the type%
\[
R_{t}=\left[
\begin{matrix}
K[z]/(z^{t}) & zK[z]/(z^{t})\\
zK[z]/(z^{t}) & K[z]/(z^{t})
\end{matrix}
\right]
\]
appeared in the description of the T-ideals $T(S)$ of $K\langle x_{1}%
,\ldots,x_{n}\rangle$ containing $T(M_{2}(K))$ (see [2]). For every unitary
algebra $S$ such that the T-ideal $T(S)$ strictly contains $T(M_{2}(K))$ there
exists a nilpotent algebra $N$ such that $T(S)=T(R_{t}\oplus N^{\natural})$
for a suitable $t$, where the algebra $N^{\natural}$ is obtained from $N$ by
formal adjoint of $1$. (Another description of $T(S)\supset T(M_{2}(K))$ was
given by Kemer [7].) The algebra $R_{3}$ appears also in noncommutative
invariant theory (see [3]). The algebra of $G$-invariants
\[
(K\langle x_{1},\ldots,x_{n}\rangle/(K\langle x_{1},\ldots,x_{n}\rangle\cap
T(S)))^{G},\quad n\geq2,
\]
is finitely generated for every finite subgroup $G$ of $GL_{n}(K)$ if and only
if $T(S)$ is not contained in $T(R_{3})$.

\noindent For a background on algebras with polynomial identity see e.g. the
book by Giambruno and Zaicev [4]. Recall that the PI-degree $\text{PIdeg}(S)$
of a PI-algebra $S$ is equal to the maximal $p$ such that the multilinear
polynomial identities of $S$ follow from the multilinear polynomial identities
of $M_{p}(K)$.

\bigskip

\noindent\textbf{6.2.Corollary.}\textit{ Let }$A$\textit{ be an }$n\times
n$\textit{ matrix over an arbitrary field }$K$\textit{ and let }$p$\textit{ be
the maximal number of equal elementary Jordan matrices in the canonical Jordan
form of }$A$\textit{ over the algebraic closure of }$K$\textit{. Then }%
\[
\text{PIdeg}(\text{Cen}(A))=p.
\]
\noindent\textbf{Proof.} Let
\[
A=\sum_{i,j=1}^{n}a_{i,j}E_{i,j}%
\]
and let $E$ be the algebraic closure of $K$. The centralizer Cen$_{K}%
(A)\subseteq M_{n}(K)$ consists of all matrices
\[
B=\sum_{i,j=1}^{n}\xi_{i,j}E_{i,j}\in M_{n}(K),\quad\xi_{i,j}\in K,
\]
such that
\[
AB=\sum_{i,k=1}^{n}\left(  \sum_{j=1}^{n}a_{i,j}\xi_{j,k}\right)  E_{i,k}%
=\sum_{i,k=1}^{n}\left(  \sum_{j=1}^{n}\xi_{i,j}a_{j,k}\right)  E_{i,k}=BA.
\]
Hence the entries $\xi_{i,j}$ of $B$ are the solutions of the system of
$n^{2}$ homogeneous linear equations
\[
\sum_{j=1}^{n}a_{i,j}\xi_{j,k}=\sum_{j=1}^{n}\xi_{i,j}a_{j,k},\quad
i,k=1,\ldots,n.
\]
The dimension $\text{dim}_{K}($Cen$_{K}(A))$ of Cen$_{K}(A)$ over $K$ is equal
to the dimension $n^{2}-r$ of the vector space of the solutions of the system,
where $r$ is the rank of the matrix of the system. Since Cen$_{E}(A)$ is
obtained from the same homogeneous linear system but considered over $E$, we
derive that
\[
\text{dim}_{E}(\text{Cen}_{E}(A))=\text{dim}_{K}(\text{Cen}_{K}(A)),
\]%
\[
\text{Cen}_{E}(A)\cong E\otimes_{K}\text{Cen}_{K}(A).
\]
It is well known that the extension of the base field preserves the PI-degree
of the algebra and we may assume that $K$ is algebraically closed. Then for a
finite dimensional $K$-algebra $S$ with Jacobson radical $J$ the PI-degree of
$S$ is equal to the maximal size of the matrix subalgebras of $S/J$. Applying
Theorem 6.1 we complete the proof. $\square$

\bigskip

\noindent The algebra $\mathcal{C}_{A}$, as presented in the proof of Theorem
6.1, has a natural $\mathbb{Z}$-grading assuming that the variable $z$ has
degree $1$. The component of degree $0$ is isomorphic to the factor algebra
$\mathcal{C}_{A}/I$. In characteristic $0$ this algebra has several remarkable
properties obtained by Giambruno and Zaicev (see [4] for detailed exposition).
It plays a key role in their result about the exponent of PI-algebras. If
$c_{n}(R)$, $n=0,1,2,\ldots$, is the codimension sequence of the PI-algebra
$R$, then
\[
\exp(R)=\lim_{n\rightarrow\infty}\sqrt[n]{c_{n}(R)}%
\]
exists and is a nonnegative integer. The algebras $\mathcal{C}_{A}/I$ are also
minimal of given exponent. If $R$ is a finitely generated PI-algebra with the
property that $\exp(R)>\exp(S)$ for any PI-algebra $S$ such that the
polynomial identities of $S$ strictly contain the polynomial identities of
$R$, then the polynomial identities of $R$ coincide with the polynomial
identities of one of the algebras $\mathcal{C}_{A}/I$.

\bigskip

\noindent Acknowledgment: The second and the third named authors wish to thank
P.N. Anh and L. Marki for fruitful consultations.

\bigskip

\noindent REFERENCES

\bigskip

\begin{enumerate}
\item Bergman, G.: \textit{Centralizers in free associative algebras}, Trans.
Amer. Math. Soc. \textbf{137} (1969), 327-344.

\item Drensky, V.:\textit{ Polynomial identities of finite dimensional
algebras,} Serdica \textbf{12} (1986), 209-216.

\item Drensky, V.:\textit{ Finite generation of invariants of finite linear
groups on relatively free algebras,} Linear and Multilinear Algebra
\textbf{35} (1993), No. 1, 1-10.

\item Giambruno, A., Zaicev, M.: \textit{Polynomial Identities and Asymptotic
Methods,} Mathematical Surveys and Monographs \textbf{122}, Amer. Math. Soc.,
Providence, RI, 2005.

\item Guralnick, R.: \textit{A note on commuting pairs of matrices}, Linear
and Multilinear Algebra \textbf{31} (1992), 71-75.

\item Guralnick, R., Sethuraman, B.A.: \textit{Commuting pairs and triples of
matrices and related varieties}, Linear Algebra Appl. \textbf{310} (2000), 139--148.

\item Kemer, A.R.: \textit{Asymptotic basis of identities of algebras with
unit of the variety }$Var(M_{2}(F))$ (Russian), Izv. Vyssh. Uchebn. Zaved.,
Mat. (1989), No. \textbf{6}, 71-76. Translation: Sov. Math. \textbf{33}
(1990), No. 6, 71-76.

\item Nelson, G.C., Ton-That, T.:\textit{ Multiplicatively closed bases for
C(A)}, Note di Matematica \textbf{26}, No. 2 (2006), 81--104.

\item Prasolov, V.V.: \textit{Problems and Theorems in Linear Algebra}, Vol.
\textbf{134} of Translation of Mathematical Monographs, American Mathematical
Society, Providence, Rhode Island, 1994.

\item Suprunenko, D.A. and Tyshkevich, R.I.:\textit{ Commutative Matrices},
Academic Press, New York and London, 1968.

\item Szigeti, J.: \textit{Linear algebra in lattices and nilpotent
endomorphisms of semisimple modules}, Journal of Algebra \textbf{319} (2008), 296--308.
\end{enumerate}

\end{document}